\documentclass[12]{article}

\usepackage{authblk,bbm,fancyhdr}
\usepackage{amsmath,amsthm,amsfonts,amstext,amssymb,amsopn,mathtools}
\usepackage{algorithm,algpseudocode,enumitem,standalone}
\usepackage{tikz}
\usetikzlibrary{shapes, arrows,decorations,automata,backgrounds,petri,positioning,arrows}
\usepackage[hidelinks]{hyperref}
\usepackage[utf8]{inputenc}

\newtheorem{theorem}{Theorem}[section]
\newtheorem{lemma}[theorem]{Lemma}
\newtheorem{proposition}[theorem]{Proposition}

\theoremstyle{remark}
\newtheorem{remark}{Remark}[section]

\usepackage[a4paper, margin=1in, headheight=13.6pt]{geometry} 

\newcommand{\RR}{\mathbb{R}}
\newcommand{\set}[1]{\mathcal{#1}}
\newcommand{\expect}[2]{\mathbb{E}_{#2}[#1]}
\newcommand{\argmax}{\mathrm{argmax}}

\newcommand{\worstQ}{\mathbb{Q}}
\newcommand{\probset}{\mathfrak{P}}

\newcommand{\rv}{\tilde{\omega}}
\newcommand{\rvset}{\Omega}
\newcommand{\obs}{\omega}

\newcommand{\M}{\mathcal{M}}

\newcommand{\distance}[1]{\text{dist}(#1)}
\newcommand{\setDistance}[1]{\mathbb{D}(#1)}
\newcommand{\hausdorffDistance}[1]{\mathbb{H}(#1)}

\newcommand{\momentApprox}{\widehat{\probset}_{\textup{mom}}}

\title{Stochastic Decomposition Method for Two-Stage Distributionally Robust Optimization}
\date{First Submitted: November 4, 2020}
\author{Harsha Gangammanavar\thanks{Department of Engineering Management, Information, and Systems, Southern Methodist University, Dallas, TX 75275 (\texttt{harsha@smu.edu}).} \text{ and} Manish Bansal\thanks{Department of Industrial and Systems Engineering, Virginia Tech, Blacksburg, VA 24061
  (\texttt{bansal@vt.edu}).}}

\begin{document}
\maketitle

{\bf Abstract.} In this paper, we present a sequential sampling-based algorithm for the two-stage distributionally robust linear programming (2-DRLP) models. The 2-DRLP models are defined over a general class of ambiguity sets with discrete or continuous probability distributions. The algorithm is a distributionally robust version of the well-known stochastic decomposition algorithm of Higle and Sen (Math. of OR 16(3), 650-669, 1991) for a two-stage stochastic linear program. We refer to the algorithm as the distributionally robust stochastic decomposition (DRSD) method. The key features of the algorithm include (1) it works with data-driven approximations of ambiguity sets that are constructed using samples of increasing size and (2) efficient construction of approximations of the worst-case expectation function that solves only two second-stage subproblems in every iteration. We identify conditions under which the ambiguity set approximations converge to the true ambiguity sets and show that the DRSD method asymptotically identifies an optimal solution, with probability one. We also computationally evaluate the performance of the DRSD method for solving distributionally robust versions of instances considered in stochastic programming literature. The numerical results corroborate the analytical behavior of the DRSD method and illustrate the computational advantage over an external sampling-based decomposition approach (distributionally robust L-shaped method).

\section{Introduction} \label{sect:intro}
Many applications of practical interest have been formulated as stochastic programming (SP) models. The models with recourse, particularly in a two-stage setting, have gained wide acceptance across application domains. These two-stage stochastic linear programs (2-SLPs) can be stated as follows:
\begin{align} \label{eq:2slp_master}
    \min~ \{c^\top x + \expect{Q(x,\rv)}{P^\star}~|~ x \in \set X\}. 
\end{align}
One needs to have complete knowledge of the probability distribution $P^\star$ to formulate the above problem. Alternatively, one must have an appropriate means to simulate observations of the random variable so that a sample average approximation (SAA) problem with a finite number of scenarios can be formulated and solved. In many practical applications, distribution associated with random parameters in the optimization model is not precisely known. It either has to be estimated from data or constructed by expert judgments, which tend to be subjective. In any case, identifying a  distribution using available information may be cumbersome at best. Stochastic min-max programming that has gained significant attention in recent years under the name of \emph{distributionally robust optimization} (DRO) is intended to alleviate the ambiguity is distributional information.%

In this paper, we study a particular manifestation of the DRO problem in the two-stage setting, viz., the two-stage distributionally robust linear programming (2-DRLP) problem. This problem is stated as:%
\begin{align}
	\min~ \{f(x) = c^\top x + \worstQ(x) ~|~ x \in \set{X} \}. \label{eq:2drlp_master}
\end{align}
Here, $c$ is the coefficient vector of a linear cost function and $\set X \subseteq \RR^{d_x}$ is the feasible set of the first-stage decision vector. The feasible region $\set{X}$ takes the form of a compact polyhedron $\set{X} = \{x ~|~ Ax \geq b, x \geq 0\}$, where $A \in \RR^{m_1 \times d_x}$ and $b \in \RR^{m_1}$. The function $\worstQ(x)$ is the worst-case expected recourse cost, which is formally defined as follows:
\begin{align}
    \worstQ(x) = \max_{P \in \probset}~ \bigg\{\mathcal{Q}(x; P) := \expect{Q(x,\rv)}{P}\bigg\}. \label{eq:distrSeparation}
\end{align}
The random vector $\rv \in \RR^{d_\obs}$ is defined on a measurable space $(\rvset, \cal F)$, where $\Omega$ is the sample space equipped with the sigma-algebra $\set{F}$ and $\probset$ is a set of continuous or discrete probability distributions defined on the measurable space $(\rvset, \set{F})$. The set of probability distributions is often referred to as the \emph{ambiguity set}. The expectation operation $\expect{\cdot}{P}$ is taken with respect to the probability distribution $P \in \probset$. For a given $x \in \mathcal{X}$, we refer to the optimization problem in \eqref{eq:distrSeparation} as the \emph{distribution separation problem}. For a given realization $\obs$ of the random vector $\rv$, the recourse cost in \eqref{eq:distrSeparation} is the optimal value of the following second-stage linear program:
\begin{align} \label{eq:2drlp_subproblem}
    Q(x,\obs) := \min \quad & g(\obs)^\top y \\
    \text{s.t.} \quad
    & y \in \set Y(x,\obs) := \big\{W(\obs) y = r(\obs) - T(\obs) x,~ y \geq 0\big\} \subset \RR^{d_y}. \notag
\end{align}
The second-stage parameters $g \in \mathbb{R}^{d_y}$, $W \in \mathbb{R}^{m \times d_y}$, $r \in \mathbb{R}^m$, and $T \in \RR^{m\times d_x}$ can depend on uncertainty.

Most data-driven approaches for 2-SLP, such as SAA, tackle the problem in two steps. In the first step, an uncertainty representation is generated using a finite set of observations that serves as an approximation of $\Omega$ and the corresponding empirical distribution serves as an approximation of $P^\star$. For a given uncertainty representation, one obtains a deterministic approximation of \eqref{eq:2slp_master}. In the second step, the approximate problem is solved using deterministic optimization methods. Such a two-step approach may lead to poor out-of-sample performance, forcing the entire process to be repeated from scratch with an improved uncertainty representation. Since sampling is performed prior to the optimization step, the two-step approach is also referred to as the \emph{external sampling procedure}.

The data-driven approaches for DRO problems avoid working with a fixed approximation of $P^\star$ in the first step. However, the ambiguity set is still defined either over the original sample space $\Omega$ or a finite approximation of it. Therefore, the resulting ambiguity set is a deterministic representation of the true ambiguity set in \eqref{eq:2drlp_master}. Once again, deterministic optimization tools are employed to solve the DRO problem. In many data-driven settings, prior knowledge of the sample space may not be available, and using a finite sample to approximate the original sample space may result in similar out-of-sample performance as in the case of the external sampling approach for 2-SLP. 

\subsection{Contributions} In light of the above observations regarding the two-step procedure for data-driven optimization, the main contributions of this manuscript are highlighted in the following.
\begin{enumerate}
    \item \emph{A Sequential Sampling Algorithm}: We present a sequential sampling approach for solving 2-DRLP. We refer to this algorithm as the \emph{distributionally robust stochastic decomposition} (DRSD) algorithm following its risk-neutral predecessor, the two-stage stochastic decomposition (SD) method \cite{Higle1991}. The algorithm uses a sequence of ambiguity sets that evolve over the course of the algorithm due to the sequential inclusion of new observations. While the simulation step improves the representation of the ambiguity set, the optimization step improves the solution in an online manner. Therefore, the DRSD method concurrently performs simulation and optimization steps. Moreover, the algorithm design does not depend on any specific ambiguity set description, and hence, is suitable for a general family of ambiguity sets.
    \item \emph{Convergence Analysis}: The DRSD method is an inexact bundle method that creates outer linearization for the dynamically evolving approximation of the first-stage problem. We provide the asymptotic analysis of DRSD and identify conditions on ambiguity sets under which the sequential sampling approach identifies an optimal solution to the 2-DRLP problem in \eqref{eq:2drlp_master} with probability one.
    \item \emph{Computational Evidence of Performance}: We provide the first evidence that illustrates the advantages of a sequential sampling approach for DRO through computational experiments conducted on well-known instances in SP literature. 
\end{enumerate}

\subsection{Related work}
The DRSD method has its roots in two-stage SP, in particular two-stage SD. In this subsection, we review the related two-stage SP literature along with decomposition and reformulation-based approaches for 2-DRLP.


For 2-SLP problems with finite support, including the SAA problem, the L-shaped method due to Van Slyke and Wets \cite{VanSlyke1969} has proven to be very effective. Other algorithms for 2-SLP's such as the Dantzig-Wolfe decomposition \cite{Dantzig1960} and the progressive hedging (PH) algorithm \cite{Rockafellar1991} also operate on problems with finite support. The well-established theory of SAA (see Chapter 5 in \cite{Shapiro2014}) supports the external sampling procedure for 2-SLP. The quality of the solution obtained by solving an SAA problem is assessed using the procedures developed, e.g., in \cite{Bayraksan2006}.  When the quality of the SAA solution is not acceptable, a new SAA is constructed with a larger number of observations. Prior work, such as \cite{Bayraksan2011} and \cite{Royset2013}, provide rules on how to choose the sequence of sample sizes in a sequential SAA procedure. 

In contrast to the above, SD-based methods incorporate one new observation in every iteration to create approximations of the dynamically updated SAA of \eqref{eq:2slp_master}. First proposed in \cite{Higle1991}, this method has seen significant development in the past three decades with the introduction of quadratic regularization term in \cite{Higle1994}, statistical optimality rules \cite{Higle1999}, and extensions to multistage stochastic linear programs \cite{Gangammanavar2020sdlp, Sen2014}. The DRSD method presented in this manuscript extends the sequential sampling approach (i.e., SD) for 2-SLPs to DRO problems. Since the simulation of new observations and optimization steps are carried out in every iteration of the SD-based methods, they can also be viewed as \emph{internal sampling methods}. 

The concept of DRO dates back to the work of Scarf \cite{Scarf1958}, and has gained significant attention in recent years. We refer the reader to \cite{rahimian2019distributionally} for a comprehensive treatment on various aspects of the DRO. The algorithmic works on DRO are either decomposition-based or reformulation-based approaches. The decomposition-based methods for 2-DRLP mimic the two-stage SP approach of using a deterministic representation of the sample space using a finite number of observations. As a consequence, the SP solution methods with suitable adaptation can be applied to solve the 2-DRLP problems. For instance, Breton and El Hachem \cite{Breton1995a, Breton1995b} apply the PH algorithm for a 2-DRLP model with a moment-based ambiguity set. Riis and Anderson \cite{Riis2005} extend the L-shaped method for 2-DRLP with continuous recourse and moment-based ambiguity set. Bansal et.al. \cite{bansal_DROdecomposition_2018} extend the algorithm in \cite{Riis2005}, which they refer to as the distributionally robust (DR) L-shaped method, to solve 2-DRLPs, with ambiguity set defined over a polytope, in finite iterations. Further extensions of this decomposition approach are presented in \cite{bansal_DROdecomposition_2018} and \cite{bansal_DRO-Disjunctive_2019} for DRO with mixed-binary recourse and disjunctive programs, respectively.

Another predominant approach to solve 2-DRLP problems is to reformulate the distribution separation problem in \eqref{eq:distrSeparation} as a minimization problem, pose the problem in \eqref{eq:2drlp_master} as a single deterministic optimization problem, and use off-the-shelf deterministic optimization tools to solve the reformulation. For example, Shapiro and Kleywegt \cite{Shapiro2002} and Shapiro and Ahmed~\cite{Shapiro2004} provided approaches for 2-DRLP problem with moment matching set to derive an equivalent stochastic program with a certain reference distribution. Bertsimas et al. \cite{Bertsimas2010} provided tight semidefinite programming reformulations for 2-DRLP where the ambiguity set is defined using multivariate distributions with known first and second moments. Likewise, Hanasusanto and Kuhn \cite{Hanasusanto2018} provided a conic programming reformulation for 2-DRLP problem where the ambiguity set comprises of a 2-Wasserstein ball centered at a discrete distribution. Xie \cite{xie2019tractable} provided similar reformulations to tractable convex programs for 2-DRLP problems with ambiguity set defined using $\infty-$ Wasserstein metric. Jiang and Guan \cite{Jiang2018} reduced the worst-case expectation in 2-DRLP, where the ambiguity set is defined using $l_1$-norm on the space of all (continuous and discrete) probability distributions, to a convex combination of CVaR and an essential supremum. By taking the dual of inner maximization, Love and Bayraksan \cite{Bayraksan2015} demonstrated that 2-DRLP where the ambiguity set is defined using $\phi$-divergence and finite sample space is equivalent to 2-SLP with a coherent risk measure. When reformulations result in equivalent stochastic programs (in \cite{Jiang2018, Bayraksan2015, ShaAhm04}, for instance), a SAA of the reformulation is used to obtain an approximate solution.

Data-driven approaches for DRO have been presented for specific ambiguity sets. In \cite{Delage2010}, problems with ellipsoidal moment-based ambiguity set whose parameters are estimated using sampled data are addressed. Esfahani et. al. tackled data-driven problems with Wasserstein metric-based ambiguity sets with convex reformulations in \cite{MohajerinEsfahani2018}. In both these works, the authors provide finite-sample performance guarantees that probabilistically bound the gap between approximate and true DRO problems. Sun and Xu presented asymptotic convergence analysis of DRO problems with ambiguity sets that are based on moments and mixture distributions constructed using a finite set of observations in \cite{Sun2016}. A practical approach to incorporate the results of these works to identify a high-quality DRO solution will be similar to the sequential SAA procedure for SP in \cite{Bayraksan2011}. Such an approach will involve the following steps performed in a series -- a deterministic representation of ambiguity set using sampled observations, applying appropriate reformulation, and solving the resulting deterministic optimization problem. If the quality of the solution is deemed insufficient, then the entire series of steps is repeated with an improved representation of the ambiguity set (possibly with a larger number of observations).

\subsubsection*{Organization} The remainder of the paper is organized as follows. In \S\ref{sect:approximations}, we present the two key ideas of the DRSD, viz., the sequential approximation of the ambiguity set and the recourse function. We provide a detailed description of the DRSD method in \S\ref{sect:algorithm}. We show the convergence of the value functions and solutions generated by the DRSD method in \S\ref{sect:convergence}. We present results of our computational experiments in \S\ref{sect:computations}, and finally we conclude and discuss about potential extensions of this paper in \S \ref{sect:conclusion}. 

\subsubsection*{Notations and Definitions}
We define the ambiguity sets over $\set{M}$, the set of all finite signed measures on the measurable space $(\rvset, \set{F})$. A nonnegative measure (written as $P \succeq 0$) that satisfies $P(\rvset) = 1$ is a probability measure. For probability distributions $P, P^\prime \in \probset$, we define \begin{align} \label{eq:zetaDistance}
    \distance{P,P^\prime} := \sup_{F \in \mathcal{F}} \Big | \expect{F(\rv)}{P} - \expect{F(\rv)}{P^\prime} \Big |
\end{align}
as the uniform distance of expectation, where $\mathcal{F}$ is a class of measurable functions. The above is the distance with $\zeta$-structure that is used for the stability analysis in SP \cite{romisch2003stability}. The distance between a single probability distribution $P$ to a set of distributions $\probset$ is given as $\distance{P,\probset} = \inf_{P^\prime \in \probset} d(P,P^\prime)$. The distance between two sets of probability distributions $\probset$ and $\widehat{\probset}$ is given as 
\begin{align*}
    \setDistance{\probset, \widehat{\probset}} := \sup_{P \in \widehat{\probset}} \distance{P, \probset}.
\end{align*}
Finally, the Hausdorff distance between $\probset$ and $\widehat{\probset}$ is defined as
\begin{align*}
    \hausdorffDistance{\probset,\widehat{\probset}} := \max\{\setDistance{\probset,\widehat{\probset}},~ \setDistance{\widehat{\probset}, \probset}\}.
\end{align*}
With suitable definitions for the set $\mathcal{F}$, the distance in \eqref{eq:zetaDistance} accepts the bounded Lipschitz, the Kantorovich and the $p$-th order Fourier-Mourier metrics (see \cite{romisch2003stability}). 



\section{Approximating Ambigiuty Set and Recourse Function} \label{sect:approximations}
In this section, we present the building blocks for the DRSD method. Specifically, we present the procedures to obtain approximations of the ambiguity set $\mathfrak{P}$ and the recourse function $Q(x,\omega)$. These procedures will be embedded within a sequential sampling-based approach. Going forward we make the following assumptions on the 2-DRLP models: 
\begin{enumerate}[label=(A\arabic{enumi})]
    \item The first-stage feasible region $\set X$ is a non-empty and compact set. \label{assum:compactX}
    \item The recourse function satisfies relatively complete recourse. The dual feasible region of the recourse problem is nonempty compact polyhedral set. The transfer (or technology) matrix satisfies $\sup_{P \in \probset} \expect{T(\rv)}{P} < \infty$.\label{assum:completeRecourse}
    \item The randomness affects the right-hand sides of constraints in \eqref{eq:2drlp_subproblem}. \label{assum:rhs}
    \item Sample space $\rvset$ is a compact metric space and the ambiguity set $\probset \neq \emptyset$.\label{assum:compactOm}
\end{enumerate}
As a consequence of \ref{assum:completeRecourse}, the recourse function satisfies $Q(x,\rv) < \infty$ with probability one for all $x \in \set{X}$. It also implies that the second-stage feasible region, i.e., $\{y~:~ W y = r(\obs) - T(\obs)x,~ y\geq 0\}$, is non-empty for all $x \in \set X$ and every $\obs \in \rvset$. The non-empty dual feasible region $\Pi$ implies that there exists a $L > -\infty$ such that $Q(x,\rv) > L$. Without loss of generality, we assume that $Q(x,\rv)\geq 0$. As a consequence of \ref{assum:rhs}, the cost coefficient vector $g$ and the recourse matrix $W$ are not affected by uncertainty. Problems that satisfy \ref{assum:compactOm} are said to have a fixed recourse. Finally, the compactness of the support $\rvset$ guarantees that every probability measure $P \in \probset$ is tight.

\subsection{Approximating the Ambiguity Set} \label{sect:ApproxAmbiguitySet}
The DRO approach assumes only partial knowledge about the underlying uncertainty that is captured by a suitable description of the ambiguity set. An ambiguity set must capture the true distribution with absolute or high degree of certainty, and must be computationally manageable. In this section we present a family of of ambiguity sets that are of interest to us in this paper.  

The computational aspects of solving a DRO problem relies heavily on the structure of the ambiguity set. The description of these structures involve parameters which are determined based on practitioner's risk-preferences. The ambiguity set descriptions that are prevalent in the literature include moment-based ambiguity sets with linear constraints (e.g., \cite{Dupacova(1987)distributionRO}) or conic constraints (e.g., \cite{Delage2010}); Kantorovich distance or Wasserstein metric-based ambiguity sets \cite{Mehrotra2013}; $\zeta$-structure metrics \cite{Zhao2015}, $\phi$-divergences such as $\chi^2$ distance and Kullback-Leibler divergence \cite{Ben-tal2012}; Prokhorov metrics \cite{Erdogan2006}, among others. Although the design of the DRSD method can work with any ambiguity set description  defined over a compact sample space, we use 2-DRLPs with moment-based ambiguity sets and Wasserstein distance-based ambiguity sets to illustrate the algorithm in details. 

In a data-driven setting, the parameters used in the description of ambiguity sets are estimated using a finite set of independent observations which can either be past realizations of the random variable $\rv$ or simulated by an oracle. We will denote such a sample by $\rvset^k \subseteq \rvset$ with $\rvset^k = \{\obs^j\}_{j=1}^k$. Naturally, we can view $\rvset^k$ as a random sample and define the empirical frequency
\begin{align}
	\hat{p}^k(\obs^j) = \frac{\kappa(\obs^j)}{k} \qquad \text{for all } \obs^j \in \rvset^k,
\end{align}
where $\kappa(\obs^j)$ denotes the number of times $\obs^j$ is observed in the sample. Since in sequential sampling setting, the sample set is updated within the optimization algorithm, it is worthwhile to note that the empirical frequency can be updated using the following recursive equations:
\begin{align} \label{eq:probUpdate}
	\hat{p}^k(\obs) = \left \{ \begin{array}{ll}
		\theta^k \hat{p}^{k-1}(\obs) & \text{if}~ \obs \in \rvset^{k-1}, \obs \neq \obs^k \\ 
		\theta^k \hat{p}^{k-1}(\obs) + (1-\theta^k) & \text{if}~ \obs \in \rvset^{k-1}, \obs = \obs^k \\
		(1-\theta^k)  & \text{if}~ \obs \notin \rvset^{k-1}, \obs = \obs^k.
	\end{array}\right. 
\end{align}
where $\theta^k \in [0,1]$. We will succinctly denote the the above operation using the mapping $\Theta^k: \RR^{|\rvset^{k-1}|} \rightarrow \RR^{|\rvset^k|}$. 

In remainder of this section, we will present alternate descriptions of ambiguity sets and show the construction of what we will refer to as {\it approximate ambiguity sets}, denoted by $\probset^k$. Let $\set{F}^k = \sigma(\obs^j~|~ j \leq k)$ be the $\sigma$-algebra generated by the observations in the sample $\rvset^k$. Notice that $\set{F}^{k-1} \subseteq \set{F}^k$, and hence, $\{\set{F}^k\}_{k \geq 1}$ is a filtration. We will define the approximate ambiguity sets over the measurable space $(\rvset^k, \set{F}^k)$. These sets should be interpreted to include all distributions that could have generated using the sample $\rvset^k$, which share a certain relationship with sample statistics. We will use $\set{M}^k$ to denote the finite signed measures on $(\rvset^k, \set{F}^k)$.

\subsubsection{Moment-based Ambiguity Sets} \label{sect:momentAmbiguity}
Given the first $q$ moments associated with the random variable $\rv$, the moment-based ambiguity set can be defined as
\begin{align} \label{eq:momentAmbuity}
	\probset_{\text{mom}} = \left \{ P \in \set{M} \left \vert \begin{array}{l}
		\int_\rvset dP(\rv) = 1, \\ \int_\rvset \psi_i(\rv)dP(\rv) = b_i \qquad i = 1,\ldots,q
	\end{array}\right. \right \}.
\end{align}
While the first constraint ensures the definition of a probability measure, the moment requirements are guaranteed by the second constraints. Here, $\psi_i(\rv)$ denote real valued measurable function on $(\rvset, \set{F})$ and $b_i \in \RR$ be a scalar for $i = 1,\ldots,q$. Existence of moments ensures that $b_i < \infty$ for all $i = 1,\ldots,q$. Notice that the description of the ambiguity set requires explicit knowledge of the following statistics: support $\rvset$ and the moments $b_i$ for $i = 1,\ldots,q$. In the data-driven setting, the support is approximated by $\rvset^k$ and the sample moments $\hat{b}_i^k = (1/k)\sum_{j=1}^k \psi_i(\obs^j)$ are used to define the following approximate ambiguity set
\begin{align} \label{eq:momentAmbuity_approx}
	\momentApprox^k = \left \{ P \in \set{M}^k \left \vert \begin{array}{l}
		\sum_{\obs \in \rvset^k} p(\obs) = 1, \\ 
		\sum_{\obs \in \rvset^k} p(\obs) \psi_i(\obs) = \hat{b}_i^k \qquad i = 1,\ldots,q
	\end{array}\right. \right \}.
\end{align}

The following result characterizes the relationship between distributions drawn from the above approximate ambiguity set, as well as asymptotic behavior of the sequence $\{\momentApprox^k\}_{k \geq 1}$.
\begin{proposition} \label{prop:momentAmbiguity_property}
For any $P \in \momentApprox^{k-1}$, we have $\Theta^k (P) \in \momentApprox^k$ where $\theta^k = \frac{k-1}{k}$. Further, suppose $\momentApprox^k \neq \emptyset$ for all $k \geq 0$, $\hausdorffDistance{\momentApprox^k, \probset_{\textup{mom}}} \rightarrow 0$ as $k \rightarrow \infty$, almost surely. 
\end{proposition}
\begin{proof}
See Appendix \S\ref{sect:proofs}.
\end{proof}

In the context of DRO, similar ambiguity sets have been studied in \cite{Bertsimas2005,Dupacova(1987)distributionRO} where only the first moment (i.e., $q = 1$) was considered. The form of ambiguity set above also relates to those used in \cite{Delage2010, Riis2005, Scarf1958, Sun2016} where constraints were imposed only on the mean and covariance. In the data driven setting of \cite{Delage2010} and \cite{Sun2016}, the statistical estimates are used in constructing the approximate ambiguity set as in the case of \eqref{eq:momentAmbuity_approx}. However, the ambiguity sets in these previous works are defined over the original sample space $\rvset$, as opposed to $\rvset^k$ that is used in \eqref{eq:momentAmbuity_approx}. This marks a critical deviation in the way the approximate ambiguity set are constructed.

\begin{remark}\label{rem:RiisandAndersonProp2.1}
When the moment information is available about the underlying distribution $P^\star$, an approximate moment-based ambiguity set with constant parameters in \eqref{eq:momentAmbuity_approx} (i.e., with $\hat{b}_i^k = b_i$ for all $k$) can be constructed. Such an approximate ambiguity sets defined over $\rvset^k$ is studied in \cite{Riis2005}. Notice that these approximate ambiguity sets satisfy $\cup_{k\geq 1} \widehat{\probset}^k \subseteq \probset$ and $\widehat{\probset}^k\subseteq \widehat{\probset}^{k+1}$, for all $k\geq 1$. 
\end{remark}

\subsubsection{Wasserstein distance-based Ambiguity Sets}  \label{sect:wassersteinAmbiguity}
We next present approximations of another class of ambiguity sets that has gained significant attention in the DRO literature, viz., the Wasserstein distance-based ambiguity sets. Consider probability distributions $\mu_1, \mu_2 \in \set{M}$, and a function $\nu:\rvset \times \rvset \rightarrow \RR_+ \cup \{\infty\}$ such that $\nu$ is symmetric, $\nu^\frac{1}{r}(\cdot)$ satisfies triangle inequality for $1 \leq r < \infty$, and $\nu(\obs_1, \obs_2) = 0$ whenever $\obs_1 = \obs_2$.  If $\Pi(\mu_1, \mu_2)$ denotes the joint distribution of random vectors $\obs_1$ and $\obs_2$ with marginals $\mu_1$ and $\mu_2$, respectively, then the Wasserstein metric of order $r$ is given by 
\begin{align} \label{eq:wassersteinMetric}
	d_{\text{w}}(\mu_1, \mu_2) = \inf_{\eta \in \Pi(\mu_1, \mu_2)} \bigg\{ \int_{\rvset \times \rvset} \nu(\obs_1, \obs_2)\eta(d\obs_1, d\obs_2) \bigg\}. 
\end{align}
In the above definition, the decision variable $\eta \in \Pi$ can be viewed as a plan to transport goods/mass from an entity whose spatial distribution is given by the measure $\mu_1$ to another entity with spatial distribution $\mu_2$. Therefore, the $d_{\text{w}}(\mu_1,\mu_2)$ measures the optimal transport cost between the measures. Notice that an arbitrary norm $\|\bullet\|^r$ on $\RR^{d_\obs}$ satisfies the requirement of the function $c(\cdot)$. In this paper, we will use $\|\bullet\|$ as the choice of our metric, in which case we obtain the Wasserstein metric of order~1. Using this metric, we define an ambiguity set as follows:
\begin{align} \label{eq:wassersteinAmbiguity}
	\probset_{\textup{w}} = \{P \in \set{M} ~|~ d_W(P,P^*) \leq \epsilon\}
\end{align}
for a given $\epsilon > 0$ and a reference distribution $P^*$. As was done in \S\ref{sect:momentAmbiguity}, we define approximate ambiguity sets defined over the measurable space $(\rvset^k, \mathcal{F}^k)$ as follows:
\begin{align} \label{eq:wassersteinAmbiguityApprox}
	\widehat{\probset}_{\textup{w}}^k = \{P \in \set{M}^k ~|~ d_W(P,\widehat{P}^k) \leq \epsilon\}.
\end{align}
For the above approximate ambiguity set, the distribution separation problem in \eqref{eq:distrSeparation} takes the following form:
\begin{subequations} \label{eq:distrSeparationWasserstein}
\begin{align}
	\max~ & \sum_{\obs \in \rvset^k} p(\obs) Q(x,\obs) \\ 
	\text{subject to}~&P \in \widehat{\probset}^k_{\textup{w}} = \left \{ P \in \set{M}^k \left \vert 
	\renewcommand{\arraystretch}{1.5}
	\begin{array}{l}
		\sum_{\obs \in \rvset^k} p(\obs) = 1 \\
		\sum_{\obs^\prime \in \rvset^k} \eta(\obs,\obs^\prime) = p(\obs) \qquad \forall \obs \in \rvset^k, \\
		\sum_{\obs \in \rvset^k} \eta(\obs,\obs^\prime) = \hat{p}^k(\obs^\prime) \qquad \forall \obs^\prime \in \rvset^k, \\
		\sum_{(\obs, \obs^\prime) \in \rvset^k \times \rvset^k} \|\obs - \obs^\prime\| \eta(\obs,\obs^\prime) \leq \epsilon \\
		\eta(\obs,\obs^\prime) \geq 0 \quad \forall \obs, \obs^\prime \in \rvset^k
	\end{array}\right. \right \}. \label{eq:wassersteinAmbiguityApprox_full}
\end{align}
\end{subequations}


The following result characterizes the distributions drawn from the approximate ambiguity sets of the form in \eqref{eq:wassersteinAmbiguityApprox}, or equivalently \eqref{eq:wassersteinAmbiguityApprox_full}.
\begin{proposition} \label{prop:wassersteinAmbiguity_property}
The sequence of Wasserstein distance-based approximate ambiguity set satisfies the following properties (1) for any $P \in \widehat{\probset}_{\textup{w}}^{k-1}$, we have $\Theta^k(P) \in  \widehat{\probset}_{\textup{w}}^{k}$ where $\theta^k = \frac{k-1}{k}$, and (2) $\mathbb{H}( \widehat{\probset}_{\textup{w}}^{k},  \probset_{\textup{w}}) \rightarrow 0$ as $k \rightarrow \infty$, almost surely.
\end{proposition}
\begin{proof}
See appendix \S\ref{sect:proofs}.
\end{proof}

The approximate ambiguity set in \cite{MohajerinEsfahani2018} is a ball constructed in the space of probability distributions that are defined over the sample space $\rvset$ and whose radius reduces with increase in the number of observations. Using Wasserstein balls of shrinking radii, the authors of \cite{MohajerinEsfahani2018} show that the optimal value of the sequence of DRO problems converges to the optimal value of the expectation-valued SP problem in \eqref{eq:2slp_master} associated with the true distribution $P^\star$.  The approximate ambiguity set in \eqref{eq:wassersteinAmbiguityApprox}, on the other hand, uses a constant radius for all $k \geq 1$. In this regard, we consider settings where the ambiguity is not necessarily resolved with increasing number of observations. This is reflected in the approximate ambiguity sets \eqref{eq:momentAmbuity_approx} and \eqref{eq:wassersteinAmbiguityApprox} that converge to their respective true ambiguity sets \eqref{eq:momentAmbuity} and \eqref{eq:wassersteinAmbiguity}, respectively.

\subsection{Approximating the Recourse Problem} \label{sect:recourseApprox}
Cutting plane methods for the 2-SLPs use an outer linearization-based approximation of the first-stage objective function in \eqref{eq:2slp_master}. In such algorithms, the challenging aspect of computing the expectation is addressed by taking advantage of the structure of the recourse problem \eqref{eq:2drlp_subproblem}. Specifically, for a given $\obs$, the recourse value $Q(\cdot,\obs)$ is known to be convex in the right-hand side parameters that includes the first-stage decision vector $x$. Additionally, if the support of $\rv$ is finite and \ref{assum:completeRecourse} holds, then the function $Q(\cdot,\obs)$ is polyhedral. Under assumptions \ref{assum:completeRecourse} and \ref{assum:compactOm}, these structural property of convexity extends to the expected recourse value $\mathcal{Q}(x)$. 

Due to strong duality of linear programs, the recourse value is also equal to the optimal value of the dual of \eqref{eq:2drlp_subproblem}, i.e., 
\begin{align} \label{eq:subproblemDual}
	Q(x,\obs) = \max~& \pi^\top [r(\obs) - T(\obs)x] \\
	\text{subject to}~& \pi \in \Pi := \{\pi~|~W^\top \pi \leq g\}. \notag 
\end{align}
Due to \ref{assum:completeRecourse} and \ref{assum:compactOm}, the dual feasible region $\Pi$ is a polytope that is not affected by the uncertainty. If $\Pi \subset \Pi$ denotes the set of all extreme points of the polytope $\Pi$, then the recourse value can also be expressed as the pointwise maximum of affine functions computed using elements of the set $\Pi$ as given below.
\begin{align} \label{eq:recoursePolyhedralForm}
	Q(x,\obs) = \max_{\pi \in \Pi} \pi^\top [r(\obs) - T(\obs)x].
\end{align}
The outer linearization approaches tend to approximate the above form of recourse function by identifying the extreme points (optimal solutions to \eqref{eq:subproblemDual}) at a sequence of candidate (or trial) solutions $\{x^k\}$, and generating the corresponding affine functions. If $\pi(x^k,\obs)$ is the optimal dual obtained by solving \eqref{eq:subproblemDual} with $x^k$ as input, then the affine function $\alpha^k(\obs) + (\beta^k(\obs))^\top x$ is obtained by computing the coefficients $\alpha^k(\obs) = (\pi^k(\obs))^\top r(\obs)$ and $\beta^k(\obs) = C(\obs)^\top \pi^k(\obs)$. Following linear programming duality, notice that this affine function is a supporting hyperplane to $Q(x,\obs)$ at $x^k$, and lower bounds the function at every other $x \in \set{X}$.

If the support $\rvset$ is finite, then one can solve a dual subproblem for all $\obs \in \rvset$ with the candidate solution as input, generate the affine functions, and collate them together to obtain an approximation of the first-stage objective function. This is the essence of the L-shaped method applied to 2-SLP in \eqref{eq:2slp_master}. In each iteration of the L-shaped method, the affine functions generated using a candidate solution $x^k$ and information gathered from individual observations are weighed by the probability density of the observation to update the approximate first-stage objective function. Notice that even when SAA of the first-stage objective function  of the 2-SLP, using a sample $\rvset_N \subset \Omega$ of size $N$, the L-shaped method can be applied. A similar approximation strategy is used in the DR L-shaped method for 2-DRLP problems. 

Alternatively, we can consider the following approximation of the recourse function expressed in the form given in \eqref{eq:recoursePolyhedralForm}:
\begin{align}\label{eq:recoursePolyhedralApprox}
	Q^k(x,\obs) = \max_{\pi \in \Pi^k} \pi^\top[r(\obs) - C(\obs) x].
\end{align}
Notice that the above approximation is built using only a subset $\Pi^k \subset \Pi$ of extreme points, and therefore, satisfies $Q^k(x,\obs) \leq Q(x,\obs)$. Since $Q(x,\obs) \geq 0$, we begin with $\Pi^0 = \{0\}$. Subsequently, we construct a sequence of sets $\{\Pi^k\}$ such that $\Pi^k \subseteq \Pi^{k+1} \subseteq \ldots \subset \Pi$ that ensures $Q^k(x,\obs)\geq 0$ for all $k$. The following result from \cite{Higle1991} captures the behavior of the sequence of approximation $\{Q^k\}$.
\begin{proposition}
The sequence $\{Q^k(x,\obs)\}_{k \geq 1}$ converges uniformly to a continuous function on $\set{X}$ for any $\obs \in \rvset$. \label{prop:uniformConvergenceApproxRecourse}
\end{proposition}
\begin{proof}
See Appendix \ref{sect:proofs}.
\end{proof}

The approximations of the form in \eqref{eq:recoursePolyhedralApprox} is one of the principal features of the SD algorithm (see \cite{Higle1991, Higle1994}). While the L-shaped and DR L-shaped methods require a finite support for $\rv$, SD is applicable even for problems with continuous support. The algorithm uses an ``incremental'' SAA for the first-stage objective function by adding one new observation in each iteration. Therefore, the first-stage objective function approximation used in SD is built using the recourse problem approximation in \eqref{eq:recoursePolyhedralApprox} and the incremental SAA. This approximation is given by:
\begin {align} \label{eq:saaIncremental}
	\mathcal{Q}^k(x) = c^\top x + \frac{1}{k} \sum_{j=1}^k Q^k(x,\obs^j). 
\end{align}
The affine functions generated in SD provide an outer linearization for the approximation in \eqref{eq:saaIncremental}. The sequence of sets that grow monotonically in size, viz. $\{\Pi^k\}$, is generated by adding one new vertex to the previous set $\Pi^{k-1}$ to obtain the updated set $\Pi^k$. The newly added vertex is the optimal dual solution obtained by solving \eqref{eq:subproblemDual} with the most recent observation $\obs^k$ and candidate solution $x^k$ as input.

We refer the reader to \cite{Birge2011}, \cite{bansal_DROdecomposition_2018,Riis2005}, and \cite{Higle1991,Higle1996} for the a detailed exposition of the L-shaped, the DR L-Shaped, and the SD methods, respectively, and note only the key differences between these methods. Firstly, the sample used to in the (DR) L-shaped method is fixed prior to optimization. In SD, this sample is updated dynamically throughout the course of the algorithm. Secondly, exact subproblem optimization for all observations in the sample is used in every iteration of the (DR) L-shaped method. On the contrary, exact optimization is used only for the subproblems corresponding to the latest observation, and an ``argmax'' procedure (to be described in the next section) is used for observations encountered in earlier iterations.

\section{Distributionally Robust Stochastic Decomposition}\label{sect:algorithm}
In this paper, we focus on a setting where the ambiguity set $\probset$ is approximated by a sequence of ambiguity sets $\{\widehat{\probset}^k\}_{k > 0}$ such that the following properties are satisfied: ($i$) for any $P \in \widehat{\probset}^{k-1}$, there exists $\theta^k \in [0,1]$ such that $\Theta^k(P) \in  \widehat{\probset}^{k}$ and ($ii$) $\mathbb{H}( \widehat{\probset}^{k},  \probset) \rightarrow 0$ as $k \rightarrow \infty$, almost surely. The moment-based ambiguity set $\probset_{\text{mom}}$ and Wasserstein distance-based ambiguity set $\probset_w$ are two sets that satisfy these  properties (Propositions \ref{prop:momentAmbiguity_property} and \ref{prop:wassersteinAmbiguity_property}, respectively). Recall that the approximate ambiguity set in an iteration (say $k$) is constructed using a finite set of observations $\rvset^k$ that progressively grow in size. Note that the sequence of approximate ambiguity sets may not necessarily converge to a single distribution. In other words, we do not assume that increasing sample size will overcome ambiguity asymptotically, as is the case in \cite{MohajerinEsfahani2018, Zhao2015}.

\begin{algorithm}[!ht]
\caption{Distributionally Robust Stochastic Decomposition}
\begin{algorithmic}[1]
    \State {\bf Input:} Incumbent solution $\hat{x}^1 \in \set{X}$; initial sample $\rvset^0 \subseteq \rvset$; stopping tolerance $\tau > 0$; $\gamma \in (0,1]$, and minimum iterations $k^{\min}$.
    \State {\bf Initialization:} Set iteration counter $k\leftarrow 1$; $\Pi^0 = \emptyset$; $\set{L}^0 = \emptyset$, and $f^0(x)  = 0$.
    \While{($k \geq k^{\min}$ and $f^{k-1}(\hat{x}^{k-1}) - f^{k-1}(x^k) \geq \tau f^{k-1}(\hat{x}^{k-1})$)}	
    \State  Solve the master problem $\M^k$ \eqref{eq:2rdsd_master} to obtain a candidate solution $x^k$. \label{step:masterProblem}
    \State Generate a scenario $\obs^k \in \rvset$ to get sample $\rvset^k \leftarrow \rvset^{k-1} \cup \{\obs^k\}$. \label{step:scenariogeneraration} 
    \For{$\obs \in \rvset^k$} \label{step:stocUpdates}
        \If {$\obs = \obs^k$} \label{step:currentObs}
            \State Solve the second-stage linear program \eqref{eq:2drlp_subproblem} with ($x^k,\obs$) as input;
            \State Obtain the optimal value $Q(x^k,\obs)$ and optimal dual solution $\pi(x^k,\obs)$;
            \State Update dual vertex set $\Pi^k \leftarrow \Pi^{k-1} \cup \{\pi(x^k,\obs)\}$. \label{step:currentObs_end}
        \Else \label{step:oldObs}
            \State Use the argmax procedure \eqref{eq:argmax} to identify dual vertex $\pi(x^k,\obs^k)$;
            \State Store $Q^k(x^k,\obs) = (\pi(x^k,\obs))^\top[r(\obs) - T(\obs) x^k]$.\label{step:oldObs_end}
        \EndIf
    \EndFor
	\State Solve the  distribution separation problem using the ambiguity set $\widehat{\probset}^k$ and $\{Q^k(x^k,\obs)\}_{\obs \in \rvset^k}$ to get an extremal distribution $P^k:= (p^k(\obs))_{\obs \in \rvset_k}$. \label{line:dissepalgo}
	\State Derive affine function ${\ell}_k^k(x) = {\alpha}_k^k + ({\beta}_k^k)^\top x$ using $\{\pi(x^k,\obs)\}_{\obs \in \rvset^{k}}$ and $P^k$ to get lower bound approximation of $\worstQ^k(x)$ as in \eqref{eq:affineCoeff};\label{step:affinefunctionlowerbound}
	\State Perform Steps \ref{step:stocUpdates}-\ref{step:affinefunctionlowerbound} with $\hat{x}^{k-1}$ (incumbent solution) to obtain $\hat{\ell}_k^k(\cdot)$. \label{step:repeatforincumbentsol}
	\For{$\ell_j^{k-1} \in \set{L}^{k-1}$} \label{step:affinecoeffupdateini}
	    \State Update previously generated affine functions $\ell^{k-1}_j(x) = {\alpha}_{k-1}^j + ({\beta}_{k-1}^j)^\top x$: $$\alpha^k_j = \theta^k \alpha^{k-1}_j \text{ and } \beta^k_j = \theta^k \alpha^{k-1}_j;$$ 
	    \State Set ${\ell}^k_j(x) = {\alpha}^k_j + ({\beta}^k_j)^\top x$ that provides lower bound approx. of $\worstQ^k(x)$;
	\EndFor \label{step:affinecoeffupdateend}
	\State Build a collection of these affine functions, denoted by $\set{L}^k$;\label{step:buildLk}
	\State Update approximation of the first-stage objective function:\label{step:update1stageapprox}
    \begin{align*}
    	c^\top x + \worstQ^k (x) \geq f^k(x) = c^\top x + \max_{j \in \mathcal{L}^k}~ \{ \alpha^k_j + (\beta^k_i)^\top x \};
    \end{align*}
    \If{Incumbent update rule \eqref{eq:incumbUpdt} is satisfied} \label{step:incumbUpdate}
	    \State $\hat{x}^k \leftarrow \hat{x}^{k-1}$ and $i_k \leftarrow i_{k-1}$;
	 \Else
    	 \State $\hat{x}^k \leftarrow x^k$ and $i_k \leftarrow k$;
	\EndIf
	\State Update the master problem by replacing $f^{k-1}(x)$ with $f^{k}(x)$ to obtain $\mathcal{M}^{k+1}$;
	\State $k \leftarrow k+1$;\label{step:incrementcounter}
	\EndWhile \label{step:incumbUpdate_end}
\State \Return {$x^k$} \label{return}
\label{Algo:return}
	\end{algorithmic}
	\label{Algo:GenTSSDP}
\end{algorithm}

The pseudocode of the DRSD method is given in Algorithm \ref{Algo:GenTSSDP}. 
We present the main steps of the DRSD method in iteration $k$ (Steps \ref{step:masterProblem}-\ref{step:incrementcounter} of Algorithm~\ref{Algo:GenTSSDP}). At the beginning of iteration $k$, we have a certain approximation of the first-stage objective function that we denote as $f^{k-1}(x)$, a finite set of observations $\rvset^{k-1}$ and an incumbent solution $\hat{x}^{k-1}$. We will use the term \emph{incumbent solution} to refer to the best solution discovered by the algorithm until iteration $k$. The solution identified in the current iteration will be referred to as the \emph{candidate solution} and denote it as $x^k$ (without $\hat{\bullet}$).

Iteration $k$ begins by first identifying the candidate solution by solving the following the master problem (Step \ref{step:masterProblem}):
\begin{align}\label{eq:2rdsd_master}
	x^k \in \arg\min~\{f^{k-1}(x) ~|~ x \in \mathcal{X}\},
\end{align}
denoted by $\mathcal{M}^k$. Following this, a new observation $\obs^k \in \rvset$ is realized, and added to the current sample of observations $\rvset^{k-1}$ to get $\rvset^k = \rvset^{k-1} \cup \{\obs^k\}$ (Step \ref{step:scenariogeneraration}). 

In order to build the first-stage objective function approximation, we rely upon the recourse function approximation presented in Section \ref{sect:recourseApprox}. For the most recent observation $\obs^k$ and the candidate solution $x^k$, we evaluate the recourse function value $Q(x^k, \obs^k)$ by solving \eqref{eq:2drlp_subproblem}, and obtain the dual optimum solution $\pi(x^k, \obs^k)$. Likewise, we obtain dual optimum solution $\pi(\hat{x}^{k-1}, \obs^k)$ by solving \eqref{eq:2drlp_subproblem} for incumbent solution $\hat{x}^{k-1}$ (Steps \ref{step:currentObs}--\ref{step:currentObs_end}).  These dual vectors are added to a set $\Pi^{k-1}$ of previously discovered optimal dual vectors. In other words, we recursively update $\Pi^k \leftarrow \Pi^{k-1}\cup \{\pi(x^k,\obs^k), \pi(\hat{x}^{k-1},\obs^k)\}$. For all other observations ($\obs \in \rvset^k,\ \obs \neq \obs^k$), we identify a dual vector in $\Pi^k$ that provides the best lower bounding approximation at $\{Q(x^k, \obs)\}$ using the following operation (Steps \ref{step:oldObs}--\ref{step:oldObs_end}):
\begin {equation} \label{eq:argmax}
\pi(x^k, \obs) \in \arg\max~ \{ \pi^{\top} [r(\obs) - T(\obs)x^k]~|~ \pi \in \Pi^k \}. 
\end {equation}
Note that the calculations in \eqref{eq:argmax} are carried out only for previous observations as $\pi(x^k, \omega^k)$ provides the best lower bound at $Q(x^k, \obs^k)$. Further, notice that $$\pi(x^k, \obs)^\top[r(\obs) - T(\obs)x^k] = Q^k(x^k,\obs),$$ the approximate recourse function value at $x^k$ defined in \eqref{eq:recoursePolyhedralApprox}, for all $\obs \in \rvset^k$, and $Q^k(x^k,\omega^k) = Q(x^k,\omega^k)$.

Using $\{Q^k(x^k,\obs^j)\}_{j=1}^k$, we solve a \emph{distribution separation problem} (in Step \ref{line:dissepalgo}): 
\begin{align} \label{eq:distrSeparationApprox}
    \worstQ^k(x^k) = \max~\bigg \{\sum_{\obs \in \rvset^k} p(\obs)Q^k(x^k,\obs)~|~ p(\obs) \in \widehat{\probset}^k \bigg \}.
\end{align}
Let $P^k = (p^k(\obs))_{\obs \in \rvset^k}$ denote the optimal solution of the above problem which we identify as the maximal/extremal probability distribution. Since the problem is solved over measures $\set{M}^k$ that are defined only over the observed set $\rvset^k$, the maximal probability distribution has weights $p^k(\obs^j)$ for $\obs^j \in \rvset^k$, and $p^k(\obs) = 0$ for $\obs \in \rvset\setminus \rvset^k$. Notice that the problem in \eqref{eq:distrSeparation} differs from the distribution separation problem \eqref{eq:distrSeparationApprox} as the latter uses the recourse function approximation $Q^k(\cdot)$ and approximate ambiguity set $\widehat{\probset}^k$ as opposed to the true recourse function $Q(\cdot)$ and ambiguity set $\probset$, respectively. For the ambiguity sets in \eqref{sect:ApproxAmbiguitySet} the distribution separation problem is a deterministic linear program. In general, the distribution separation problems associated with well-known ambiguity sets remain deterministic convex optimization problems \cite{rahimian2019distributionally}, and off-the-shelf solvers can used to obtain the extremal distribution.


In Step \ref{step:affinefunctionlowerbound} of Algorithm \ref{Algo:GenTSSDP}, we use the dual vectors $\{\pi(x^k,\obs^j)\}_{j \leq k}$ and the maximal probability distribution $P^k$ to generate a lower bounding affine function:
\begin{align} \label{eq:cutComputation}
    \worstQ^k(x) = \max_{P \in \widehat{\probset}^k} \expect{Q^k(x,\rv)}{P} \geq \sum_{\obs^j \in \rvset^k} p^k(\obs^j) \cdot (\pi(x^k,\obs^j))^{\top} [r(\obs^j) - T(\obs^j)x],
\end{align}
for the worst case expected recourse function measured with respect to the maximal probability distribution $P^k \in \widehat{\probset}^k$ which is obtained by solving the distribution separation problem \eqref{eq:distrSeparationApprox}. We denote the coefficients of the affine function on the right-hand side of \eqref{eq:cutComputation} by
\begin{align}\label{eq:affineCoeff}
	\alpha_k^k = \sum_{\obs^j \in \rvset^k} p^k(\obs^j) \pi(x^k,\obs^j)^\top r(\obs^j)\text{ and }\beta_k^k = -\sum_{\obs^j \in \rvset^k} p^k(\obs^j) T(\obs^j)^\top \pi(x^k,\obs^j),
\end{align}
 and succinctly write the affine function as $\ell_k^k(x) = \alpha_k^k + (\beta_k^k)^\top x$. Similar calculations are carried out using the incumbent solution $\hat{x}$ to identify a maximal probability distribution and a lower bounding affine function  resulting in the affine function $\hat{\ell}_k^k(x) = \hat{\alpha}_k^k + (\hat{\beta}_k^k)^\top x$. 

While the latest affine functions provide lower bound on $\worstQ^k$, the affine functions generated in previous iteration are not guaranteed to lower bound $\worstQ^k$. To see this, let us consider the moment-based approximate ambiguity sets $\{\widehat{\probset}^k_{\text{mom}}\}_{k\geq 1}$. Let $P^j_{\text{mom}} \in \widehat{\probset}^j_{\text{mom}}$ be the maximal distribution identified in an iteration $j <k$ which was used to compute the affine function $\ell_j^j(x)$. By assigning $p^j(\obs) = 0$ for all new observations encountered after iteration $j$, i.e., $\obs \in \rvset^k \setminus \rvset^j$, we can construct a probability distribution $\bar{P} = ((p^j(\obs))_{\obs \in \rvset^j}, (0)_{\obs \in \rvset^k \setminus \rvset^j}) \in \mathbb{R}^{|\Omega^k|}_+$. This reconstructed distribution satisfies $\sum_{\obs \in \rvset^k} \bar{p}(\obs) = 1$. However, it is easy to see that it does not satisfy $\sum_{\obs \in \rvset^k} \psi_i(\obs) \bar{p}(\obs) = \hat{b}^j_i = \hat{b}^k_i$ for all $i = 1,\ldots,q$. Therefore, $\bar{P} \notin \probset^k$. In other words, while the coefficients $(\alpha_j^j, \beta_j^j)$ are $\set{F}^j$-measurable, the corresponding measure is not feasible to the approximate ambiguity set $\probset^k$. Therefore, $\ell_j^j(x)$ is not a valid lower bound to $\worstQ^k$. The arguments for the Wasserstein-based approximate ambiguity set are more involved, but persistence of a similar issue can be demonstrated.

In order to address this, we update the previously generated affine functions $\ell^{k-1}_j(x) = {\alpha}_{k-1}^j + ({\beta}_{k-1}^j)^\top x$ for $j<k$, as follows (Steps \ref{step:affinecoeffupdateini} - \ref{step:affinecoeffupdateend}):
\begin{align} \label{eq:affineCoeff_update}
	\alpha_j^k = \theta^k \alpha_j^{k-1}, \ \ \beta_j^k = \theta^k \beta_j^{k-1}, \text{ and } {\ell}^k_j(x) = {\alpha}_k^j + ({\beta}_k^j)^\top x \qquad  \text{ for all } j < k,
\end{align}
such that $\ell^k_j(x)$ provides lower bound approximation of $\worstQ^k(x)$ for all $j \in \{1,\ldots, k-1\}$. Similarly, we update the affine functions $\hat{\ell}^k_j(x)$, $j<k$, associated with incumbent solution (Step \ref{step:repeatforincumbentsol}). Recall that for $\widehat{\probset}_{\text{mom}}$ and $\widehat{\probset}_w$, the parameter $\theta^k = \frac{k-1}{k}$ (Propositions \ref{prop:momentAmbiguity_property} and \ref{prop:wassersteinAmbiguity_property}). The candidate and the incumbent affine functions ($\ell_k^k(x)$ and $\hat{\ell}_k^k(x)$, respectively), as well as the updated collection of previously generated affine functions are used to build the set of affine functions which we will denote by $\set{L}^k$ (Step \ref{step:buildLk}). The lower bounding property of this first-stage objective function approximation is captured in the following result.
\begin{theorem} \label{thm:lowerBoundingMinorants}
Under assumption \ref{assum:completeRecourse}, the first-stage objective function approximation in \eqref{eq:objfnApprox} satisfies 
\begin{align*}
f^k(x) \leq c^\top x + \worstQ^k(x) \text{ for all } x \in \set{X} \text{ and } k \geq 1.
\end{align*}
\end{theorem}
\begin{proof}
For non-empty approximate ambiguity set $\widehat{\probset}^1$ of ambiguity set $\probset$, the construction of the affine function ensures that $\ell_1^1(x) \leq \worstQ^1(x)$. Assume that $\ell(x) \leq \worstQ^{k-1}(x)$ for all $\ell \in \set{L}^{k-1}$ and $k>1$. The maximal nature of the probability distribution $P^k$ satisfies:
\begin{align*}
	\sum_{\obs \in \rvset^k} p^k(\obs) Q^k(x,\obs) &\geq \sum_{\obs \in \rvset^k} p(\obs) Q^k(x,\obs) \qquad \forall P \in \widehat{\probset}^k.
\end{align*}
Using above and the monotone property of the approximate recourse function, we have 
\begin{align} 
	\sum_{\obs \in \rvset^k} p^k(\obs) Q^k(x,\obs) &\geq \sum_{\obs \in \rvset^k} p(\obs) Q^{k-1}(x,\obs) \notag \\
	&= \sum_{\obs \in \rvset^k\setminus \obs^k} p(\obs) Q^{k-1}(x,\obs) + p(\obs^k) Q^{k-1}(x,\obs^k), \label{eq:thm_monotonicity}
\end{align}
for all $\{p(\omega)\}_{\omega \in \Omega^k} \in \widehat{\probset}^k$. Based on the properties of $\probset$ and $\{\widehat{\probset}^k\}_{k \geq 1}$ (similar to Propositions \ref{prop:momentAmbiguity_property} and \ref{prop:wassersteinAmbiguity_property}), we know that for every $P \in \widehat{\probset}^{k-1}$ we can construct a probability distribution in $\probset^k$ using the mapping $\Theta^k$ defined by \eqref{eq:probUpdate}. Considering a probability distribution $P' = \{p'(\omega)\}_{\omega \in \Omega^{k-1}} \in \widehat{\probset}^{k-1}$ such that $\Theta^k(P') \in \widehat{\probset}^k$, the inequality \eqref{eq:thm_monotonicity} reduces to
%
\begin{align*}
	\sum_{\obs \in \rvset^k} p^k(\obs) Q^k(x,\obs) &\geq \hspace{-0.4em} \sum_{\obs \in \rvset^k\setminus \obs^k} [\theta^k p^\prime(\obs) Q^{k-1}(x,\obs)] + [\theta^k p^\prime(\obs^k) + (1-\theta^k)] Q^{k-1}(x,\obs^k) \\
	&= \theta^k \bigg [\sum_{\obs \in \rvset^{k-1}} p^\prime(\obs) Q^{k-1}(x,\obs)\bigg] + (1-\theta^k) Q^{k-1}(x,\obs^k) \\
	&\geq \theta^k \bigg [\sum_{\obs \in \rvset^{k-1}} p^\prime(\obs) Q^{k-1}(x,\obs)\bigg].
\end{align*}
The last inequality is due to assumption $\ref{assum:completeRecourse}$, i.e., $Q(x,\omega^k) \geq 0$ and the construction of recourse function approximation $Q^k$ described in \S\ref{sect:recourseApprox}. Since $\ell(x)$ lower bounds the term in bracket, we have
\begin{align*}
	\sum_{\obs \in \rvset^k} p^k(\obs) Q^k(x,\obs) &\geq \theta^k \ell(x).
\end{align*}
Using the same arguments for all $\ell \in \set{L}^{k-1}$, and the fact that  the $\ell_k^k(x)$ and $\hat{\ell}(x)$ are constructed as lower bounds to the $\worstQ^k$, we have $f^k(x) \leq c^\top x + \worstQ^k(x)$. This completes the proof by induction.
\end{proof}

Using the collection of affine functions $\set{L}^k$, we update approximation of the first-stage objective function in Step~\ref{step:update1stageapprox}, as follows:
\begin{align}\label{eq:objfnApprox}
	f^k(x) = c^\top x + \max_{i \in \mathcal{L}^k}~ \{ \alpha^i + (\beta^i)^\top x \}.
\end{align}
Once the approximation is updated, the performance of the candidate solution is compared relative to the incumbent solution (Steps \ref{step:incumbUpdate}--\ref{step:incumbUpdate_end}). This comparison is performed by verifying if the following inequality
\begin{align} \label{eq:incumbUpdt}
    f^k(x^k) - f^k(\hat{x}^{k-1}) < \gamma[f^{k-1}(x^k) - f^{k-1}(\hat{x}^{k-1}],
\end{align}
where parameter $\gamma \in (0,1]$, is satisfied. If so, the candidate solution is designated to be the next incumbent solution, i.e., $\hat{x}^{k} = x^k$. If the inequality is not satisfied, the precious incumbent solution is retained as $\hat{x}^k = \hat{x}^{k-1}$. This completes an iteration of the DRSD method.

\begin{remark}
The algorithm design can be extended to incorporate 2-DRLP where the relatively complete recourse assumption of \ref{assum:completeRecourse} and/or assumption \ref{assum:rhs} is not satisfied. For problems where relatively complete recourse condition is not met, a candidate solution may lead to one or more subproblems to be infeasible. In this case, the dual extreme rays can be used to compute a feasibility cut that is included in the first-stage approximation. The argmax procedure in \eqref{eq:argmax} is only valid when assumption \ref{assum:rhs} is satisfied. In problems where the uncertainty also affects the cost coefficients, the argmax procedure presented in \cite{gangammanavar2020stochastic} can be utilized. These algorithmic enhancements can be incorporated without affecting the convergence properties of DRSD that we present in the next section. 
\end{remark}


\section{Convergence Analysis} \label{sect:convergence}
In this section we provide the convergence result of the sequential sampling-based approach to solve DRO problems. In order to facilitate the exposition of our theoretical results, we will define certain quantities for notational convenience that are not necessarily computed during the course of the algorithm in the form presented in the previous section. Our convergence results are built upon stability analyses presented in \cite{Sun2016} and convergence analysis of the SD algorithm in \cite{Higle1991}. 

We define a function which is defined over the approximate ambiguity set using the recourse function $Q(\cdot, \cdot)$, that is
\begin{align} \label{eq:objfnApprox_trueRecourse}
	g^k(x) := c^\top x + \max_{P \in \widehat{\probset}^k} \expect{Q(x,\rv)}{P}
\end{align}
for a fixed $x \in \set{X}$. We begin by analyzing the behavior of the sequence $\{g^k\}_{k \geq 1}$ as $k \rightarrow \infty$. In particular, we will assess the sequence of function evaluations at a converging subsequence of first-stage solutions. The result is captured in the following proposition.

\sloppy
\begin{proposition}\label{prop:FnConvergenceApproxSet}
Suppose $\{\hat{x}^{k_n}\}$ denotes a subsequence of $\{\hat{x}^k\}$ such that $\hat{x}^{k_n} \rightarrow \bar{x}$, then $\lim_{n \rightarrow \infty} |g^{k_n}(\hat{x}^{k_n}) - f(\bar{x})| = 0$, with probability one.
\end{proposition}
\begin{proof}
Consider the ambiguity set $\widehat{\probset}^k$. For $i = 1,2$ and $x_i \in \mathcal{X}$, let $P(x_i) \in \argmax_{P \in \widehat{\probset}^k} \{\expect{Q(x_i,\rv)}{P}$. Then, 
\begin{align*}
	g^k(x_1) =~& c^\top x_1 + \expect{Q(x_1,\rv)}{P(x_1)} \\
	\geq~& c^\top x_1 + \expect{Q(x_1,\rv)}{P(x_2)} \\
	=~& c^\top x_2 + \expect{Q(x_2,\rv}{P(x_2)} + c^\top (x_1 - x_2) + \\  & \hspace{3cm} \expect{Q(x_1,\rv)}{P(x_2)} - \expect{Q(x_2,\rv)}{P(x_2)} \\
	=~& g^k(x_2) + c^\top (x_1 - x_2) + \expect{Q(x_1,\rv)}{P(x_2)} - \expect{Q(x_2,\rv)}{P(x_2)}.
\end{align*}
The inequality in the above follows from optimality of $P(x_1)$. The above implies that 
\begin{align}
g^k(x_2) - g^k(x_1) \leq~& c^\top (x_2-x_1) +  \expect{Q(x_2,\rv)}{P(x_2)} - \expect{Q(x_1,\rv)}{P(x_2)} \notag \\
	 \leq~& |c^\top (x_2-x_1)| + \bigg | \expect{Q(x_2,\rv)}{P(x_2)} - \expect{Q(x_1,\rv)}{P(x_2)} \bigg | \notag \\
	\leq~& (\|c\| + C)\|x_2 - x_1\|. \label{eq:continuity_g1}
\end{align}
The second relationship is due to the triangular inequality. The third inequality follows from uniform Lipschitz continuity of recourse function $Q(x,\rv)$, with probability one, under assumption \ref{assum:completeRecourse} which implies that there exists a constant $C$ such that $| \expect{Q(x_1,\rv)}{P} - \expect{Q(x_2,\rv)}{P} | \leq C\|x_1 - x_2\|$ for any $P$. Therefore, the function $g^k(x)$ is equi-continuous on $x \in \set{X}$.  Starting with $x_2$ and using the same arguments, we have
\begin{align}
	g^k(x_1) - g^k(x_2) \leq~& (\|c\| + C)\|x_2 - x_1\|. \label{eq:continuity_g2}
\end{align}
Now consider ambiguity sets $\probset$ and $\widehat{\probset}^k$. Note that
\begin{align*}
	|f(x) - g^k(x)| =&~ \bigg| \max_{P \in \probset} \expect{Q(x,\rv)}{P} - \max_{P^\prime \in \widehat{\probset}^k} \expect{Q(x,\rv)}{P^\prime} \bigg| \qquad \forall x \in \set{X}\\
	\leq &~ \max_{P \in \probset} \min_{P^\prime \in \widehat{\probset}^k} \big| \expect{Q(x,\rv)}{P} - \expect{Q(x,\rv}{P^\prime}\big| \qquad \forall x \in \set{X} \\ 
	\leq&~ \max_{P \in \probset} \min_{P^\prime \in \widehat{\probset}^k} \sup_{x \in \set{X}} \big| \expect{Q(x,\rv)}{P} - \expect{Q(x,\rv)}{P^\prime}\big|.
\end{align*}
Using the definition of deviation between ambiguity sets $\probset$ and $\widehat{\probset}^k$ as well as its relationship with Hausdorff distance, we have
\begin{align}
	|f(x) - g^k(x)| = \mathbb{D}(\probset, \widehat{\probset}^k)  \leq \mathbb{H}(\probset, \widehat{\probset}^k) \label{eq:g_vs_f}.
\end{align}
For $\hat{x}^{k_n}$ and $\bar{x}$, combining \eqref{eq:continuity_g1}, \eqref{eq:continuity_g2}, and \eqref{eq:g_vs_f}, we have
\begin{align*}
	|f(\bar{x}) - g^{k_n}(\hat{x}^{k_n})| \leq&~ |f(\bar{x}) - g^{k_n}(\bar{x})| + |g^{k_n}(\bar{x}) - g^{k_n}(\hat{x}^{k_n})| \\
	\leq&~ \mathbb{H}(\probset, \widehat{\probset}^k) + (\|c\| + C) \|\bar{x} - \hat{x}^{k_n}\|.
\end{align*}
As $n \rightarrow \infty$, the family of ambiguity sets considered satisfy $\mathbb{H}(\probset, \widehat{\probset}^{k_n}) \rightarrow 0$ and the hypothesis informs us that $\hat{x}^{k_n} \rightarrow \bar{x}$. Therefore, we conclude that $g^{k_n}(\hat{x}^{k_n}) \rightarrow f(\bar{x})$ as $n \rightarrow \infty$.
\end{proof}

Notice that, the behavior of the approximate ambiguity sets defined in \S\ref{sect:ApproxAmbiguitySet}, in particular, the condition $\mathbb{H}(\probset, \widehat{\probset}^k) \rightarrow 0$ as $k \rightarrow \infty$ plays a central role in the above proof. Recall that for the moment and Wasserstein distance-based ambiguity sets, the condition is established in propositions \ref{prop:momentAmbiguity_property} and \ref{prop:wassersteinAmbiguity_property}, respectively. It is also worthwhile to note that under the foregoing conditions, \eqref{eq:g_vs_f} also implies uniform convergence of the sequence $\{g^k\}$ to $f(x)$, with probability one.

The above result applies to any algorithm that generates a converging sequence of iterates $\{x^k\}$ and a corresponding sequence of extremal distributions. Such an algorithm is guaranteed to exhibit convergence to the optimal distributionally robust objective function value. Therefore, this result is applicable to the sequence of instances constructed using external sampling and solved, for example, using reformulation-based methods. Such an approach was adopted in \cite{Riis2005} and \cite{Sun2016}. The analysis in \cite{Riis2005} relies upon two rather restrictive assumptions. The first assumption is that for all $P \in \probset$ there exists a sequence of measures $\{P^k\}$ such that $P^k \in \widehat{\probset}^k$ and converges weakly to $P$. The second assumption requires the approximate ambiguity sets to be strict subsets of the true ambiguity set, i.e., $\widehat{\probset}^k \subset \probset$. Both of these assumptions are very difficult to satisfy in a data-driven setting (also see Remark \ref{rem:RiisandAndersonProp2.1}). 

The analysis in \cite{Sun2016}, on the other hand, does not make the above assumptions. Therefore, their analysis is more broadly applicable in settings where external sampling is used to generate $\rvset^k$. DRO instances are constructed based on statistics estimated using $\rvset^k$ and solved to optimality for each $k \geq 1$. They show the convergence of optimal objective function values and optimal solution sets of approximate problems to the optimal objective function value and solutions of the true DRO problem, respectively. In this regard, the result in Proposition \ref{prop:FnConvergenceApproxSet} can alternatively be derived using Theorem 1(i) in \cite{Sun2016}. While the above function is not computed during the course of the sequential sampling algorithm, it provides the necessary benchmark for our convergence analysis.

One of the main point of deviation in our analysis stems from the fact that we use the objective function approximations that are built based on the approximate recourse function in \eqref{eq:recoursePolyhedralApprox}. In order to study the piecewise affine approximation of the first-stage objective function, we introduce another benchmark function
\begin{align}\label{eq:objfnApprox_approxRecourse}
	\phi^k(x) := c^\top x + \max_{P \in \widehat{\probset}^k} \expect{Q^k(x,\rv)}{P}. 
\end{align}
Notice that the above function uses the approximations for the ambiguity set (as in the case of \eqref{eq:objfnApprox_trueRecourse}) as well as the approximation for recourse function. This construction ensures for all $x \in \set{X}$ and $k \geq 1$ that $\phi^k(x) \leq g^k(x)$ that follows from the fact that $Q^k(x,\rv) \leq Q(x,\rv)$, almost surely. Further, following the result in Theorem \ref{thm:lowerBoundingMinorants} ensures that $f^k(x) \leq \phi^k(x)$. Putting these together, we obtain the following relationship
\begin{align}
	f^k(x) \leq \phi^k(x) \leq g^k(x) \qquad \forall x \in \set{X}, k \geq 1.
\end{align}
While the previous proposition was focused on the upper limit in the above relationship, we present the asymptotic behavior of the $\{f^k\}$ sequence in the following results.

\begin{lemma}\label{lemma:asymptoticSupport}
Suppose $\{\hat{x}^{k_n}\}$ denotes a subsequence of $\{\hat{x}^k\}$ such that $\hat{x}^{k_n} \rightarrow \bar{x}$, $\lim_{n \rightarrow \infty} f^{k_n}(\hat{x}^{k_n}) - f(\bar{x}) = 0$, with probability one.
\end{lemma}
\begin{proof}
From Proposition \ref{prop:FnConvergenceApproxSet}, we have $\lim_{n \rightarrow \infty} |f(\bar{x}) - g^{k_n}(\hat{x}^{k_n})| = 0$. Therefore, there exists $N_1 < \infty$ and $\epsilon_1>0$ such that 
\begin{align}\label{eq:support_1}
\bigg | \max_{P \in \probset} \expect{Q(\bar{x},\rv)}{P} - \max_{P \in \widehat{\probset}^{k_n}} \expect{Q(\hat{x}^{k_n},\rv)}{P} \bigg| < \epsilon_1/2 \quad \forall n > N_1.
\end{align}
Now consider,
\begin{align*}
	\max_{P \in \widehat{\probset}^{k_n}} \expect{Q(\hat{x}^{k_n}, \rv)}{P} &- \max_{P \in \widehat{\probset}^{k_n}} \expect{Q^{k_n}(\hat{x}^{k_n},\rv)}{P}  \\ 
	= & \max_{P \in \widehat{\probset}^{k_n}} \big( \expect{Q(\hat{x}^{k_n}, \rv)}{P} - \expect{Q^{k_n}(\hat{x}^{k_n},\rv)}{P}\big) \\
	\leq &~ \max_{P \in \widehat{\probset}^{k_n}} \big | \expect{Q(\hat{x}^{k_n},\rv)}{P} - \expect{Q^{k_n}(\hat{x}^{k_n},\rv)}{P} \big | \\
	=&~\max_{P \in \widehat{\probset}^{k_n}} \expect{|Q(\hat{x}^{k_n},\rv) - Q^{k_n}(\hat{x}^{k_n},\rv)|}{P}.
\end{align*}
The last equality follows from the fact that $Q(x,\rv) \geq Q^k(x,\rv)$ for all $x \in \set{X}$ and $k \geq 1$, almost surely. Moreover, because of the uniform convergence of $\{Q^k\}$ (Proposition \ref{prop:uniformConvergenceApproxRecourse}), the sequence of approximate functions $\{\phi^k\}$ also convergences uniformly. This implies that, there exists $N_2 < \infty$ such that 
\begin{align} \label{eq:support_2}
	\bigg | \max_{P \in \widehat{\probset}^{k_n}} \expect{Q(\hat{x}^{k_n}, \rv)}{P} - \max_{P \in \widehat{\probset}^{k_n}} \expect{Q^{k_n}(\hat{x}^{k_n},\rv)}{P} \bigg | < \epsilon_1/2.
\end{align}
Let $N = \max\{N_1, N_2\}$. Using \eqref{eq:support_1} and \eqref{eq:support_2}, we have for all $n > N$
\begin{align*}
	\bigg| \max_{P \in \probset} \expect{Q(\bar{x},\rv)}{P} - \max_{P \in \widehat{\probset}^{k_n}} \expect{Q^{k_n}(x^{k_n},\rv)}{P} \bigg| < \epsilon_1.
\end{align*}
This implies that $|f(\bar{x}) - \phi^{k_n}(\hat{x}^{k_n})| \rightarrow 0$ as $n \rightarrow \infty$. Based on \eqref{eq:argmax}, we have $Q^{k_n}(\hat{x}^{k_n}, \obs) = (\pi(\hat{x}^{k_n},\obs))^\top [r(\obs) - T(\obs)\hat{x}^{k_n}] \geq (\pi(\hat{x}^{k_n}, \obs))^\top [r(\obs) - T(\obs) x]$ for all $x \in \set{X}$ and $\obs \in \rvset^{k_n}$. Let
\begin{align*}
	\alpha^{k_n}_{k_n} = \sum_{\obs \in \rvset^{k_n}} p^{k_n}(\obs) (\pi(\hat{x}^{k_n},\obs))^\top r(\obs)\text{ and }\beta^{k_n}_{k_n} = -\sum_{\obs \in \rvset^{k_n}} p^{k_n}(\obs) T(\obs)^\top \pi(\hat{x}^{k_n},\obs),
\end{align*}
where $\{p^{k_n}(\obs)\}_{\obs \in \rvset^{k_n}}$ is an optimal solution of the distributional separation problem \eqref{eq:distrSeparationApprox} where index $k$ is replaced by $k_n$. Then, the affine function $\alpha^{k_n}_{k_n} + (c+\beta^{k_n}_{k_n})^\top x$ provides a lower bound approximation for function $\phi^{k_n}(x)$, i.e., 
\begin{align*}
	\phi^{k_n}(x) \geq \alpha^{k_n}_{k_n} + (c+\beta^{k_n}_{k_n})^\top x \qquad \text{ for all } x \in \set{X},
\end{align*}
with strict equality holding only at $\hat{x}^{k_n}$. Therefore, using the definition of $f^k(x)$ we have $\lim_{n \rightarrow \infty} \alpha^{k_n}_{k_n} + (c+\beta^{k_n}_{k_n})^\top \hat{x}^{k_n} = \lim_{n \rightarrow \infty} f^{k_n}(\hat{x}^{k_n}) = \lim_{n \rightarrow \infty} \phi^{k_n}(\hat{x}^{k_n})  = f(\bar{x})$, almost surely. This completes the proof.
\end{proof}

The above result characterizes the behavior of the sequence of affine functions generated during the course of the algorithm. In particular, the sequence $\{f^k(\hat{x}^k)\}_{k \geq 1}$ accumulates at the objective value of the original DRO problem \eqref{eq:2drlp_master}. Recall that the candidate solution $x^k$ is a minimizer of $f^{k-1}(x)$ and an affine function is generated at this point such that $f^k(x^k) = \phi^k(x^k)$  in all iterations $k \geq 1$. However, due to the update procedure in \eqref{eq:affineCoeff_update} the quality of the estimates at $x^k$ gradually diminishes leading to a large value for $(\phi^k(x^k) - f^k(x^k))$ as $k$ increases. This emphasizes the role of the incumbent solution and computing the incumbent affine function $\hat{\ell}(x)$ during the course of the algorithm. By updating the incumbent solution and frequently reevaluating the affine functions at the incumbent solution, we can ensure that the approximation is ``sufficiently good'' in the neighborhood of the incumbent solution. In order to assess the improvement of approximation quality, we define
\begin{align}\label{eq:estError}
	\delta^k := f^{k-1}(x^k) - f^{k-1}(\hat{x}^{k-1}) \leq 0 \qquad \forall k \geq 1.
\end{align}
The inequality follows from the optimality of $x^k$ with respect to the objective function $f^{k-1}$. The quantity $\delta^k$ measures the error in objective function estimate at the candidate solution with respect to the estimate at the current incumbent solution. The following result captures the asymptotic behavior of this error term.

\begin{lemma}\label{lemma:vanishingError}
Let $\mathcal{K}$ denotes a sequence of iterations where the incumbent solution changes. There exists a subsequence of iterations, denoted as $\set{K}^* \subseteq \mathcal{K}$, such that $\lim_{k \in \set{K}^*} \delta^k = 0$.
\end{lemma} 
\begin{proof} We will consider two cases depending on whether the set $\set{K}$ is finite or not. First, suppose that $|\set{K}|$ is not finite. By the incumbent update rule and \eqref{eq:estError},
\begin{align*}
    f^{k_n}(x^{k_n}) - f^{k_n}(\hat{x}^{k_n-1}) < \gamma [f^{k_n-1}(x^{k_n}) - f^{{k_n}-1}(\hat{x}^{k_n-1})] = \gamma \delta^{k_n} \leq 0 \qquad \forall k_n \in \set{K}.
\end{align*}
Subsequently, we have $\limsup_{n \rightarrow \infty} \delta^{k_n} \leq 0$. Since $x^{k_n} = \hat{x}^{k_n}$ and $\hat{x}^{k_n-1} = \hat{x}^{k_{n-1}}$, we have
\begin{align*}
    f^{k_n}(\hat{x}^{k_n}) - f^{k_n}(\hat{x}^{k_{n-1}}) \leq \gamma \delta^{k_n} \leq 0. 
\end{align*}
The left-hand side of the above inequality captures the improvement in the objective function value at the current incumbent solution over the previous incumbent solution. Using the above, we can write the average improvement attained over $n$ incumbent changes as\sloppy
\begin{align*}
    \frac{1}{n} \sum_{j = 1}^n \bigg[f^{k_j}(\hat{x}^{k_j}) - f^{k_j}(\hat{x}^{k_{j-1}}) \bigg] \leq \frac{1}{n} \sum_{j = 1}^n \gamma \delta^{k_j} \leq 0 \qquad \text{ for all } n.
\end{align*}
This implies that \sloppy
\begin{align*}
    \frac{1}{n}\underbrace{\bigg(f^{k_n}(\hat{x}^{k_n}) - f^{k_1}(\hat{x}^{k_0}) \bigg)}_{(a)}+ \frac{1}{n}\bigg[\sum_{j=1}^{n-1} \underbrace{\bigg(f^{k_j}(\hat{x}^{k_j}) - f^{k_{j+1}}(\hat{x}^{k_j}) \bigg)}_{(b)} \bigg]  \leq \frac{1}{n} \sum_{j = 1}^n \gamma \delta^{k_j} \leq 0, \ \ \forall n.
\end{align*}
Under the assumption that the dual feasible region is non-empty and bounded (this is ensured by relatively complete recourse, \ref{assum:completeRecourse}), $\{f^k\}$ is a sequence of Lipschitz continuous functions. This along with compactness of $\set{X}$ \ref{assum:compactX}, implies that $f^{k_n}(\hat{x}^{k_n}) - f^{k_1}(\hat{x}^{k_0})$ is bounded from above. Hence, the term (a) reduces to zero as $n \rightarrow \infty$. The term (b) converges to zero, with probability one, due to uniform convergence of $\{f^k\}$. Since $\gamma \in (0,1]$, we have
\begin{align*}
    \lim_{n \rightarrow \infty} \frac{1}{n} \sum_{j = 1}^n \delta^{k_j} = 0
\end{align*}
with probability one. Further,
\begin{align*}
    \lim_{n \rightarrow \infty} \frac{1}{n} \sum_{j = 1}^n \delta^{k_j} \leq \limsup_{n \rightarrow \infty} \delta^{k_n} \leq 0.
\end{align*}
Thus, there exists a subsequence indexed by the set $\set{K}^*$ such that $\lim_{k \in \set{K}^*} \delta^k = 0$, with probability one.  

Now if $|\set{K}|$ is finite, then there exists $\hat{x}$ and $K < \infty$ such that for all $k \geq K$, we have $\hat{x}^k = \hat{x}$. Notice that, if $\lim_{k \in \set{K}^*} x^k = \bar{x}$, uniform convergence of the sequence $\{f^k\}$ and Lemma \ref{lemma:asymptoticSupport} ensure that 
\begin{subequations} \label{eq:limitApprox}\begin{align}
\lim_{k \in \set{K}^*} f^k(x^k) = \lim_{k \in \set{K}^*} f^{k-1}(x^k) = f(\bar{x}) \\ 
\lim_{k \in \set{K}^*} f^k(\hat{x}) = \lim_{k \in \set{K}^*} f^{k-1}(\hat{x}) = f(\hat{x}). 
\end{align} \end{subequations}
Further, since the incumbent is not updated in iterations $k \geq K$, we must have from the update rule in \eqref{eq:incumbUpdt} that  
\begin{align*}
	f^{k}(x^k) - f^{k}(\hat{x}) \geq \gamma [f^{k-1}(x^k) - f^{k-1}(\hat{x})] = \gamma \delta^k \quad \text{ for all } k \geq K.
\end{align*}
Using \eqref{eq:limitApprox}, we have
\begin{align*}
	\lim_{k \in \set{K}^*} \big( f^{k}(x^k) - f^{k}(\hat{x})\big) &\geq \gamma \lim_{k \in \set{K}^*} \big(f^{k-1}(x^k) - f^{k-1}(\hat{x})\big), 
\end{align*}
which implies
\begin{align*}
 f(\bar{x}) - f(\hat{x}) &\geq \gamma(f(\bar{x}) - f(\hat{x})).
\end{align*}
Since $\gamma \in (0,1]$, we must have $f(\bar{x}) - f(\hat{x}) = 0$. Hence, $\lim_{k \in \set{K}^*} \delta^k = f(\bar{x}) - f(\hat{x}) = 0$, with probability one.
\end{proof}

Equipped with the results in lemmas \ref{lemma:asymptoticSupport} and \ref{lemma:vanishingError}, we state the main theorem which establishes the existence of a subsequence of incumbent solution sequence for which every accumulation point is an optimal solution to \eqref{eq:2drlp_master}. 
\begin{theorem}
Let $\{x^k\}_{k=1}^\infty$ and $\{\hat{x}^k\}_{k=1}^\infty$ be the sequence candidate and incumbent solutions generated by the algorithm. There exists a subsequence $\{\hat{x}^k\}_{k \in \set{K}}$ for which every accumulation point is an optimal solution of 2-DRLP \eqref{eq:2drlp_master}, with probability one.
\end{theorem}
\begin{proof}
Let $x^* \in \set{X}$ be an optimal solution of \eqref{eq:2drlp_master}. Consider a subsequence indexed by $\set{K}$ such that $\lim_{k \in \set{K}} \hat{x}^k = \bar{x}$. Compactness of $\set{X}$ ensures the existence of accumulation point $\bar{x} \in \set{X}$ and therefore,
\begin{align}  \label{eq:lbMinorant_1}
	f(x^*) \leq f(\bar{x}).
\end{align}
From Theorem \ref{thm:lowerBoundingMinorants}, we have
\begin{align*}
	f^k(x) \leq~& c^\top x + \worstQ^k(x) \\
	\leq~& c^\top x + \max_{P \in \widehat{\probset}^k} \expect{Q(x,\rv)}{P} = g^k(x) \qquad \forall k, x \in \set{X}.
\end{align*}
Thus, using the uniform convergence of $\{g^k\}$ (Proposition \ref{prop:FnConvergenceApproxSet}) we have
\begin{align} \label{eq:lbMinorant_2}
	\limsup_{k \in \set{K}^\prime} f^k(x^*) \leq \lim_{k \in \set{K}^\prime} g^k(x^*) = f(x^*)
\end{align}
for all subsequences indexed by $\set{K}^\prime \subseteq\{1,2,\ldots\}$, with probability one. Recall that,
\begin{align*}
    \delta^{k} = f^{{k}-1}(x^{k}) - f^{{k}-1}(\hat{x}^{{k}-1})
    &\leq f^{{k}-1}(x^*) - f^{{k}-1}(\hat{x}^{{k}-1}) \qquad \text{for all } k \geq 1.
\end{align*}
The inequality in the above follows from the optimality of $x^{k}$ with respect to $f^{k-1}(x)$. Taking limit over $\set{K}$, we have 
\begin{align*}
	\lim_{k \in \set{K}} \delta^{k} \leq~& \lim_{k \in \set{K}}  \big(f^{k-1}(x^*) - f^{{k}-1}(\hat{x}^{{k}-1})\big) \\
	\leq~& \limsup_{k \in \set{K}} f^{k-1}(x^*) - \liminf_{k \in \set{K}} f^{{k}-1}(\hat{x}^{{k}-1}) \\
	\leq~& f(x^*) - f(\bar{x}).
\end{align*}
The last inequality follows from \eqref{eq:lbMinorant_2} and $\lim_{k \in \set{K}} f^{k-1}(\hat{x}^{k-1}) = f(\bar{x})$ (Lemma \ref{lemma:asymptoticSupport}). From Lemma \ref{lemma:vanishingError}, there exists a subsequence indexed by $\set{K}^* \subseteq \mathcal{K}$ such that $\lim_{k \in \set{K}^*} \delta^{k} = 0$. Therefore, if $\{\hat{x}^k\}_{k \in \set{K}^*} \rightarrow \bar{x}$, we have 
\begin{align*}
	f(x^*) - f(\bar{x}) \geq 0.
\end{align*}
Using \eqref{eq:lbMinorant_1} and the above inequality, we conclude that $\bar{x}$ is an optimal solution with probability one. 
\end{proof}

\section{Computational Experiment} \label{sect:computations}
In this section, we evaluate the effectiveness and efficiency of the DRSD method in solving 2-DRLP problems. For our preliminary experiments, we consider 2-DRLP problems with moment-based ambiguity set $\mathfrak{P}_{\text{mom}}$ for the first two moments ($q=2$). We report results from the computational experiments conducted on four well-known SP test problems: {\tt CEP}, {\tt PGP}, {\tt BAA}, and {\tt STORM}. The test problems have supports of size $216$, $576$, $10^{18}$ and $10^{81}$, respectively.

We use an external sampling-based approach as a benchmark for comparison. The external sampling-based instances involve constructing approximate problems of the form \eqref{eq:objfnApprox_trueRecourse} with a pre-determined number of observations $N \in \{100, 250, 500, 1000, 2000\}$ (that might not be unique). The resulting instances are solved using the DR L-Shaped method. Specifically, we compare the solution quality provided by these methods along with the solution time. 

We conduct $30$ independent replications for each problem instance with different seeds for the random number generator. The algorithms are implemented in the C programming language, and the experiments are conducted on a $64$-bit Intel core i7 - 4770 CPU at $3.4$GHz $\times ~ 8$ machine with $32$ GB memory. All linear programs, i.e., master problem, subproblems, and distribution separation problem, are solved using CPLEX 12.10 callable subroutines. The results from the experiments are presented in Tables \ref{tab:estimates} and \ref{tab:times}. The results for the instances with a finite sample size obtained from the DR L-shaped method are labeled as DRLS-$N$, where $N$ denotes the number of observations used to approximate the ambiguity set. The DRSD method is run for a minimum of $256$ iteration, while no such limit is imposed on the DR L-shaped method.

\begin{table}[!t]\centering \renewcommand{\arraystretch}{1}
    \resizebox{0.95\columnwidth}{!}{%
    \begin{tabular}{|c|c|c|c|}
        \hline
        Method & \# Iterations & Objective Estimate & \# Unique Obs. \\ 
        \hline 
        \multicolumn{4}{|c|}{{\tt PGP}} \\ \hline
        DRLS-100 & $17.867$ ($\pm 0.92$) & $457.610$ ($\pm 3.28$) & $38.233$ ($\pm 1.10$) \\
        DRLS-250 & $20.467$ ($\pm 0.73$) & $462.922$ ($\pm 2.28$) & $53.267$ ($\pm 1.74$) \\
        DRLS-500 & $20.467$ ($\pm 0.63$) & $464.704$ ($\pm 1.95$) & $68.667$ ($\pm 1.62$) \\
        DRLS-1000 & $20.500$ ($\pm 0.84$) & $466.104$ ($\pm 1.78$) & $85.833$ ($\pm 1.78$) \\
        DRLS-2000 & $20.900$ ($\pm 0.57$) & $469.579$ ($\pm 2.53$) & $102.867$ ($\pm 1.87$) \\
        \hline 
         DRSD & $503.667$ ($\pm 687.01$) & $463.185$ ($\pm 16.28$) & $65.667$ ($\pm 10.04$) \\
        \hline 
        \multicolumn{4}{|c|}{{\tt CEP}} \\ \hline
        DRLS-100 & $2.600$ ($\pm 0.19$) & $658817.129$ ($\pm 14457.30$) & $80.700$ ($\pm 1.21$) \\
        DRLS-250 & $2.267$ ($\pm 0.17$) & $680735.606$ ($\pm 10511.47$) & $147.900$ ($\pm 2.05$) \\
        DRLS-500 & $2.067$ ($\pm 0.09$) & $683252.019$ ($\pm 5948.62$) & $195.167$ ($\pm 1.44$) \\
        DRLS-1000 & $2.000$ ($\pm 0.00$) & $679665.728$ ($\pm 4926.88$) & $214.100$ ($\pm 0.47$) \\
        DRLS-2000 & $2.000$ ($\pm 0.00$) & $680744.118$ ($\pm 3872.39$) & $215.967$ ($\pm 0.07$) \\
        \hline 
         DRSD & $257.000$ ($\pm 0.00$) & $681772.359$ ($\pm 10348.32$) & $150.033$ ($\pm 2.12$) \\
        \hline 
        \multicolumn{4}{|c|}{{\tt BAA}} \\ \hline
        DRLS-100 & $247.000$ ($\pm 5.17$) & $248677.974$ ($\pm 1238.19$) & $100.000$ ($\pm 0.00$) \\
        DRLS-250 & $240.233$ ($\pm 5.74$) & $249173.795$ ($\pm 770.38$) & $250.000$ ($\pm 0.00$) \\
        DRLS-500 & $236.067$ ($\pm 6.28$) & $249827.472$ ($\pm 499.59$) & $500.000$ ($\pm 0.00$) \\
        DRLS-1000 & $229.067$ ($\pm 6.03$) & $250640.029$ ($\pm 355.53$) & $1000.000$ ($\pm 0.00$) \\
        DRLS-2000 & $219.900$ ($\pm 5.36$) & $251405.806$ ($\pm 243.54$) & $2000.000$ ($\pm 0.00$) \\
        \hline 
         DRSD & $316.367$ ($\pm 31.53$) & $250235.142$ ($\pm 737.21$) & $315.367$ ($\pm 31.53$) \\
        \hline 
        \multicolumn{4}{|c|}{{\tt STORM}} \\ \hline
        DRLS-100 & $11.667$ ($\pm 0.51$) & $15742456.082$ ($\pm 12191.85$) & $100.000$ ($\pm 0.00$) \\
        DRLS-250 & $11.167$ ($\pm 0.52$) & $15781724.510$ ($\pm 8753.59$) & $250.000$ ($\pm 0.00$) \\
        DRLS-500 & $11.733$ ($\pm 0.59$) & $15797019.902$ ($\pm 5345.83$) & $500.000$ ($\pm 0.00$) \\
        DRLS-1000 & $11.767$ ($\pm 0.52$) & $15806575.387$ ($\pm 3771.91$) & $1000.000$ ($\pm 0.00$) \\
        DRLS-2000 & $11.900$ ($\pm 0.46$) & $15817039.917$ ($\pm 2502.88$) & $2000.000$ ($\pm 0.00$) \\
        \hline 
         DRSD & $516.200$ ($\pm 108.53$) & $15786864.662$ ($\pm 9155.49$) & $515.200$ ($\pm 108.53$) \\
        \hline 
    \end{tabular}%
    }
    \caption{Comparison of results obtained from DR L-Shaped method and DRSD method}
    \label{tab:estimates}
\end{table}

Table \ref{tab:estimates} shows the average number of iterations, the average objective function value, and the average number of unique observations in Columns 2, 3, and 4, respectively. The values in the parenthesis are the half-widths of the corresponding confidence intervals. Notice that for DRSD, the number of iterations is also equal to the number of observations used to approximate the ambiguity set. The results show the ability of DRSD to dynamically determine the number of observations by assessing the progress made during the algorithm. The objective function estimate obtained using the DRSD is comparable to the objective function estimate obtained using the DR L-shaped method of comparable size. For instance, the DRSD objective function estimate for {\tt STORM} that is based upon a sample of size $516.2$ (on average) is within $0.1\%$ of the objective function value estimate of DRLS-500. These results show that the optimal objective function estimate obtained from DRSD are comparable to those obtained using an external sampling-based approach.

\begin{table}[!t]
    \centering \renewcommand{\arraystretch}{1}
    \resizebox{0.95\textwidth}{!}{%
    \begin{tabular}{|c|c|c|c|c|c|c|}
        \hline
        Method & Total & Master & Subproblem & Optimality & Argmax & Separation \\ 
         \hline 
        \multicolumn{7}{|c|}{{\tt PGP}} \\ \hline
        DRLS-100 & $0.0517$ & $0.0019$ & $0.0422$ & $0.0000$ & $-$ & $0.0023$ \\ 
        DRLS-250 & $0.0770$ & $0.0021$ & $0.0651$ & $0.0000$ & $-$ & $0.0025$ \\ 
        DRLS-500 & $0.0959$ & $0.0020$ & $0.0822$ & $0.0000$ & $-$ & $0.0026$ \\ 
        DRLS-1000 & $0.1214$ & $0.0021$ & $0.1050$ & $0.0000$ & $-$ & $0.0030$ \\ 
        DRLS-2000 & $0.1463$ & $0.0022$ & $0.1270$ & $0.0000$ & $-$ & $0.0033$ \\ 
        \hline  
         DRSD & $0.3498$ & $0.1225$ & $0.0549$ & $0.0003$ & $0.0018$ & $0.0697$ \\ 
        \hline 
        \multicolumn{7}{|c|}{{\tt CEP}} \\ \hline
        DRLS-100 & $0.0154$ & $0.0002$ & $0.0086$ & $0.0000$ & $-$ & $0.0004$ \\ 
        DRLS-250 & $0.0237$ & $0.0001$ & $0.0124$ & $0.0000$ & $-$ & $0.0004$ \\ 
        DRLS-500 & $0.0283$ & $0.0001$ & $0.0138$ & $0.0000$ & $-$ & $0.0003$ \\ 
        DRLS-1000 & $0.0285$ & $0.0001$ & $0.0134$ & $0.0000$ & $-$ & $0.0003$ \\ 
        DRLS-2000 & $0.0286$ & $0.0001$ & $0.0137$ & $0.0000$ & $-$ & $0.0003$ \\ 
        \hline 
         DRSD & $0.1154$ & $0.0299$ & $0.0262$ & $0.0001$ & $0.0009$ & $0.0292$ \\ 
        \hline 
        \multicolumn{7}{|c|}{{\tt BAA}} \\ \hline
        DRLS-100 & $2.6497$ & $0.0778$ & $2.1698$ & $0.0067$ & $-$ & $0.1837$ \\ 
        DRLS-250 & $6.1376$ & $0.0756$ & $5.2157$ & $0.0078$ & $-$ & $0.3252$ \\ 
        DRLS-500 & $11.1848$ & $0.0655$ & $9.7257$ & $0.0056$ & $-$ & $0.4624$ \\ 
        DRLS-1000 & $21.1755$ & $0.0616$ & $18.5514$ & $0.0043$ & $-$ & $0.8286$ \\ 
        DRLS-2000 & $44.2579$ & $0.0685$ & $39.2027$ & $0.0056$ & $-$ & $1.2710$ \\ 
        \hline 
         DRSD & $2.0943$ & $0.1570$ & $0.0577$ & $0.0002$ & $0.0550$ & $0.9003$ \\ 
        \hline 
        \multicolumn{7}{|c|}{\tt STORM} \\ \hline
        DRLS-100 & $0.4336$ & $0.0022$ & $0.3197$ & $0.0001$ & $-$ & $0.0391$ \\ 
        DRLS-250 & $1.0080$ & $0.0021$ & $0.7429$ & $0.0001$ & $-$ & $0.0885$ \\ 
        DRLS-500 & $2.1167$ & $0.0025$ & $1.5629$ & $0.0001$ & $-$ & $0.1957$ \\ 
        DRLS-1000 & $4.3179$ & $0.0027$ & $3.1438$ & $0.0001$ & $-$ & $0.4571$ \\ 
        DRLS-2000 & $9.0015$ & $0.0027$ & $6.2653$ & $0.0001$ & $-$ & $1.3203$ \\ 
        \hline 
         DRSD & $30.4394$ & $0.7815$ & $0.3084$ & $0.0003$ & $0.4781$ & $23.8544$ \\ 
        \hline 
        \end{tabular}%
        }
    \caption{Computational time comparison between DR L-shaped and DRSD}
    \label{tab:times}
\end{table}

Table \ref{tab:times} shows the average total computational time (Column 2) for each instance. The table also includes the average time spent to solve the Master problem (first-stage approximate problem) and the subproblems, to verify the optimality conditions, to complete the argmax procedure (only for DRSD), and to solve the distribution separation problem (Columns 3--7, respectively). The results for small scale instances ({\tt PGP} and {\tt CEP}) show that both DRSD and the DR L-shaped method take a fraction of a second, but the computational time for DRSD is higher than the DR L-shaped method for all $N$. This behavior can be attributed to the fact that (i) the subproblems are relatively easy to solve and the computational effort to solve all the subproblems does not increase significantly with $N$, and (ii) the DRSD is run for a minimum number of iterations (256) and thereby, contributing to the total time taken to solve master problems and distribution separation problems. This observation is in-line with our computational experience with the SD method for 2-SLPs. It is important to note that, while the computational time for the DR L-shaped method on an individual instance may be lower, but the iterative procedure necessary to identify a sufficient sample size may require solving several instances with an increasing sample size. This may result in a significantly higher cumulative computational time. The DRSD method, and the sequential sampling idea in general, mitigates the need for this iterative process. 

On the other hand, for large-scale problems ({\tt BAA} and {\tt STORM}), we observe a noticeable increase in the computational time for the DR L-shaped method with an increase in $N$. A significant portion of this time is spent on solving the subproblems. Since the DRSD solves only two subproblems in each iteration, the time taken to solve the subproblems is significantly less in comparison to the DR L-shaped method for all $N$ where all subproblems corresponding to unique observations are solved in each iteration. Notice that for {\tt STORM}, the average number of iterations taken by DRSD is at least $43$ times the average number of iterations taken by DR L-shaped for each $N$. This is reflected in the significant increase in the computational time for solving the master problems and the distributional separation problems. In contrast, for {\tt BAA}, the average number of DRSD iterations is only $28\%$ higher than the DR L-shaped iterations. As a result, the increase in the computational times due to an increase in the master and distributional separation problem solution times is overshadowed by the computational gains attained by solving only two subproblems in each iteration. Moreover, the computational time associated with solving the distribution separation problem can be reduced by using column-generation procedures that take advantage of the problem structure. Such an implementation was not undertaken for our current experiments and is a fruitful future research avenue.

\section{Conclusions} \label{sect:conclusion}

We presented a new decomposition approach for solving two-stage distributionally robust linear programs (2-DRLPs) with a general ambiguity set that is defined using continuous and/or discrete probability distributions with a very large sample space. Since this approach extended the stochastic decomposition approach of Higle and Sen \cite{Higle1991} for 2-DRLPs with a singleton ambiguity set, we referred to it as Distributionally Robust Stochastic Decomposition (DRSD) method. The DRSD is a sequential sampling-based approach that allowed sampling within the optimization step where only two second-stage subproblems are solved in each iteration. In contrast, an external sampling procedure utilizes the distributionally robust L-shaped method~\cite{bansal_DROdecomposition_2018} for solving 2-DRLP with a finite number of scenarios, where in each iteration all subproblems are solved. While the design of DRSD accommodates general ambiguity sets, we provided its asymptotic convergence analysis for a family of ambiguity sets that includes the well-known moment-based and Wasserstein metric-based ambiguity sets. Furthermore, we performed computational experiments to evaluate the efficiency and effectiveness of solving distributionally robust variants of four well-known stochastic programming test problems that have supports of size ranging from $216$ to $10^{81}$. Based on our results, we observed that the objective function estimates obtained using the DRSD and the DR L-shaped method are statistically comparable. These DRSD estimates are obtained while providing computational improvements that are critical for large-scale problem instances.

\appendix
\section{Proofs} \label{sect:proofs}
In this appendix, we provide the proofs for the propositions related to the asymptotic behavior of the ambiguity sets approximations defined in \S\ref{sect:ApproxAmbiguitySet} and the recourse function approximation presented in \S\ref{sect:recourseApprox}.
\begin{proof} (Proposition \ref{prop:momentAmbiguity_property})
For $P \in \momentApprox^{k-1}$, it is easy to verify that $p^\prime = \Theta P$ satisfies the support constraint, viz., $\sum_{\obs \in \rvset^k} p^\prime(\obs) = 1$. Now consider for $i = 1,\ldots,q$, we have 
\begin{align*}
	\sum_{\obs \in \rvset^{k}} p^\prime(\obs) \psi_i(\obs) =~& \sum_{\obs \in \rvset^{k-1}, \obs \neq \obs^k} p^\prime(\obs) \psi_i(\obs) + p^\prime(\obs^k) \psi_i(\obs^k) \\
	=~& \theta^k \sum_{\obs \in \rvset^{k-1}, \obs \neq \obs^k} p(\obs) \psi_i(\obs) + \theta^kp(\obs^k)\psi(\obs^k) + (1-\theta^k)\psi_i(\obs^k) \\
	=~& \theta^k \sum_{\obs \in \rvset^{k-1}} p(\obs) \psi_i(\obs) + (1-\theta^k)\psi_i(\obs^k) \\
	=~& \hat{b}_i^{k-1} + (1-\theta^k)\psi_i(\obs^k) \\
	=~& \sum_{\obs \in \rvset^{k-1}} \hat{p}^{k-1}(\obs) \psi_i(\obs) + (1-\theta^k)\psi_i(\obs^k) = \hat{b}^k_i.
\end{align*}
This implies that $\Theta^k(P) \in \momentApprox^k$.

Using Proposition 4 in \cite{Sun2016}, there exists a positive constant $\chi$ such that 
\begin{align*}
0 \leq	\hausdorffDistance{\momentApprox^k, \probset_{\textup{mom}}} \leq \chi \|\hat{\mathbf{b}}^k - \mathbf{b}\|.
\end{align*}
Here, $\mathbf{b} = (b_i)_{i=1}^q$ and $\hat{\mathbf{b}}^k = (\hat{b}_i^k)_{i=1}^q$, and $\|\cdot\|$ denotes the Euclidean norm. Since the approximate ambiguity sets are constructed using independent and identically distributed samples of $\rv$, using law of large numbers, we have $\hat{b}_i^k \rightarrow b_i$ for all $i = 1,\ldots,q$. This completes the proof.
\end{proof}

\begin{proof} (Proposition \ref{prop:wassersteinAmbiguity_property})
Consider approximate ambiguity sets $\widehat{\probset}^{k-1}_\text{w}$ and $\widehat{\probset}^k_\text{w}$ of the form given in \eqref{eq:wassersteinAmbiguityApprox_full}. Let $P^\prime = (p^\prime(\obs))_{\obs \in \Omega^{k-1}} \in \widehat{\probset}^{k-1}_\text{w}$, and let the reconstructed probability distribution be denoted by $P$. We can easily check that $P = \Theta(P^\prime)$ is indeed a probability distribution. With $P$ fixed, it suffices now to show that the polyhedron
\begin{align} \label{eq:etaPolyhedron}
\set{E}(P, \widehat{P}^k) = \left \{ \eta \in \RR^{\rvset^k \times \rvset^k} \left \vert 
	\renewcommand{\arraystretch}{1.5}
	\begin{array}{l}		
		\sum_{\obs^\prime \in \rvset^k} \eta(\obs,\obs^\prime) = p(\obs) \qquad \forall \obs \in \rvset^k, \\
		\sum_{\obs \in \rvset^k} \eta(\obs,\obs^\prime) = \hat{p}^k(\obs^\prime) \qquad \forall \obs^\prime \in \rvset^k, \\
		\sum_{(\obs, \obs^\prime) \in \rvset^k \times \rvset^k} \|\obs - \obs^\prime\| \eta(\obs,\obs^\prime) \leq \epsilon
		\end{array}\right. \right \}.
\end{align}
is non-empty. Since $P^\prime \in \widehat{\probset}^{k-1}_\text{w}$, there exist $\eta^\prime(\obs, \obs^\prime)$ for all $(\obs, \obs^\prime) \in \rvset^{k-1} \times \rvset^{k-1}$ such that the constraints in the description of the approximate ambiguity set in \eqref{eq:wassersteinAmbiguityApprox_full} are satisfied. We will do this by analyzing two possibilities, 
\begin{enumerate}
	\item We encounter a previously seen observation, i.e., $\obs^k \in \rvset^{k-1}$ and $\rvset^k = \rvset^{k-1}$. Let $\eta(\obs, \obs^\prime) = \theta^k \eta^\prime(\obs, \obs^\prime)$ for $\obs,\obs^\prime \in \rvset^{k-1}$ and $\obs \neq \obs^\prime \neq \obs^k$; and $\eta(\obs^k, \obs^k) = \theta^k \eta^\prime(\obs^k, \obs^k) + (1-\theta^k)$. We will verify the feasibility of this choice by verifying the three sets of constraints in \eqref{eq:etaPolyhedron}.
	\begin{align*}
		\sum_{\obs^\prime \in \rvset^k} \eta(\obs,\obs^\prime) =~& \sum_{\obs^\prime \in \rvset^k \setminus \{\obs^k\}} \eta(\obs, \obs^\prime) + \eta(\obs, \obs^k)\\
		=~& \sum_{\obs^\prime \in \rvset^{k-1} \setminus \{\obs^k\}} \theta^k \eta^\prime(\obs, \obs^\prime) + \theta^k \eta^\prime(\obs, \obs^k) + \mathbf{1}_{\obs = \obs^k} (1-\theta^k) \\
		=~& \theta^k \bigg ( \sum_{\obs^\prime \in \rvset^{k-1}} \eta^\prime(\obs, \obs^\prime)\bigg) + \mathbf{1}_{\obs = \obs^k} (1-\theta^k) \\
		=~& \theta^k p^\prime (\obs) + \mathbf{1}_{\obs = \obs^k} (1-\theta^k) = p(\obs) \qquad \forall \obs \in \rvset^k.
	\end{align*}
	\begin{align*}
		\sum_{\obs \in \rvset^k} \eta(\obs,\obs^\prime) =~& \sum_{\obs \in \rvset^k \setminus \{\obs^k\}} \eta(\obs,\obs^\prime) + \eta(\obs^k,\obs^\prime) \\
		=~& \sum_{\obs \in \rvset^{k-1} \setminus \{\obs^k\}} \theta^k \eta^\prime(\obs,\obs^\prime) + \theta^k \eta(\obs^k,\obs) + \mathbf{1}_{\obs^\prime = \obs^k} (1-\theta^k) \\
		=~& \theta^k \bigg(\sum_{\obs \in \rvset^{k-1}} \eta(\obs,\obs^\prime)\bigg) + \mathbf{1}_{\obs^\prime = \obs^k} (1-\theta^k) \\
		=~& \theta^k \hat{p}^{k-1}(\obs^\prime) + 	\mathbf{1}_{\obs^\prime = \obs^k} (1-\theta^k) = \hat{p}^k(\obs^\prime) \qquad \forall \obs^\prime \in \rvset^k.	
	\end{align*}
	Finally, 
	\begin{align*}
	\sum_{(\obs, \obs^\prime) \in \rvset^k \times \rvset^k} \|\obs - \obs^\prime\| & \eta(\obs,\obs^\prime) \\ 
	=~& \sum_{\substack{(\obs, \obs^\prime) \in \rvset^{k-1} \times \rvset^{k-1}\\\obs \neq \obs^\prime \neq \obs^k}} \theta^k \|\obs - \obs^\prime\| \eta(\obs,\obs^\prime) + \|\obs^k - \obs^k\| \eta(\obs^k,\obs^k) \\
	=~& \theta^k \bigg ( \sum_{(\obs, \obs^\prime) \in \rvset^{k-1} \times \rvset^{k-1}} \|\obs - \obs^\prime\| \eta(\obs,\obs^\prime) \bigg ) \leq \theta^k \epsilon \leq \epsilon.
	\end{align*}
	Since all the three constraints are satisfied, the chosen values for $\eta$ is an element of the polyhedron $\set{E}$, and therefore, $\set{E} \neq \emptyset$. 
	\item We encounter a new observation, i.e., $\obs^k \notin \rvset^{k-1}$. Let $\eta(\obs,\obs^\prime) = \theta^k \eta^\prime(\obs, \obs^\prime)$ for $\obs, \obs^\prime \in \rvset^{k-1}$, $\eta(\obs^k,\obs^\prime) = 0 $ for $\obs^\prime \in \Omega^{k-1}$, $ \eta(\obs,\obs^k) = 0$ for $\obs \in \Omega^{k-1}$, and $\eta(\obs^k,\obs^k) = (1-\theta^k)$. Let us verify the three conditions defining \eqref{eq:etaPolyhedron} with this choice for $\eta$ variables.
	\begin{align*}
	\sum_{\obs^\prime \in \rvset^k} \eta(\obs,\obs^\prime) =~& \sum_{\obs^\prime \in \rvset^k\setminus \{\obs^k\}} \eta(\obs,\obs^\prime) + \eta(\obs,\obs^k) \\	
	=~& \sum_{\obs^\prime \in \rvset^{k-1}} \theta^k \eta^\prime(\obs, \obs^\prime) + \mathbf{1}_{\obs = \obs^k} (1-\theta^k) \\ 
	=~& \theta^k p^\prime(\obs) + \mathbf{1}_{\obs = \obs^k} (1-\theta^k) = p(\obs). \\
	\sum_{\obs \in \rvset^k} \eta(\obs,\obs^\prime) =~& \sum_{\obs \in \rvset^k \setminus \{\obs^k\}} \eta(\obs,\obs^\prime) + \eta(\obs^k,\obs^\prime) \\
	=~& \sum_{\obs^\prime \in \rvset^{k-1}} \theta^k \eta^\prime(\obs, \obs^\prime) + \mathbf{1}_{\obs^\prime = \obs^k} (1-\theta^k) \\ 
	=~& \theta^k \hat{p}^{k-1} +  + \mathbf{1}_{\obs^\prime = \obs^k} (1-\theta^k) = \hat{p}^k(\obs^\prime)
	\end{align*}
Consider,
\begin{align*}
	\sum_{(\obs, \obs^\prime) \in \rvset^k \times \rvset^k} \|\obs - \obs^\prime\| & \eta(\obs,\obs^\prime) \\ =~& \sum_{(\obs, \obs^\prime) \in \rvset^{k-1} \times \rvset^{k-1}} \theta^k \|\obs - \obs^\prime\| \eta^\prime(\obs,\obs^\prime) + \|\obs^k - \obs^k\| \eta(\obs^k,\obs^k) \\
	&\hspace{1cm} + \sum_{\obs \in \Omega^k}\|\obs - \obs^k\| \eta(\obs,\obs^k) + \sum_{\obs^\prime \in \Omega^k}\|\obs^k - \obs^\prime\| \eta(\obs^k,\obs^\prime) \\
	\leq~& \theta^k \epsilon \leq \epsilon. 
\end{align*}
Therefore, the value of the $\eta$ variables satisfy the constraints and $\set{E} \neq \emptyset$. This implies that $\Theta^k(P) \in \widehat{\probset}_{\textup{w}}^k$. 
\end{enumerate}

Next, let us consider a distribution $Q \in \widehat{\probset}_{\textup{w}}^{k}$. Then, 
\begin{align*}
	d_{\textup{w}}(Q,P^*) \leq~& d_{\textup{w}} (Q,\widehat{P}^k) + d_{\textup{w}}(\widehat{P}^k,P^*) \leq \epsilon + d_{\textup{w}}(\widehat{P}^k,P^*).
\end{align*}
The above inequality is a consequence of the triangle inequality of Wasserstein distance. Since $P \in \widehat{\probset}^k_{\textup{w}}$, we have $d_{\textup{w}}(P^*,P) \leq \epsilon$. 
Under compactness assumption for $\rvset$, $d > 2$, and $\expect{\text{exp}(\|\rv\|^a)}{P^*} < \infty$, Theorem 2 in \cite{Fournier2015rate} guarantees
\begin{align*}
	\text{Prob}\big[d_{\textup{w}}(\widehat{P}^k,P^*) \leq \delta\big]  \leq \left\{
	\begin{array}{ll} C~ \text{exp}(-ck\delta^d) & \text{if}~ \delta > 1 \\
	C~ \text{exp}(-ck\delta^a) & \text{if}~\delta \leq 1
	\end{array} \right .
\end{align*} 
for all $k \geq 1$. This implies that the $\lim_{k \rightarrow \infty} d_{\textup{w}}(\widehat{P}^k,P^*) = 0$, almost surely. Consequently, we obtain that $d_{\textup{w}}(Q,P^*) \leq \epsilon$ (or equivalently $Q \in \probset_{\textup{w}}$) as $k \rightarrow \infty$, almost surely. This completes the proof.
\end{proof}

\begin{proof}(Proposition \ref{prop:uniformConvergenceApproxRecourse})
Recall that $\set{X} \times \rvset$ is a compact set because of Assumptions \ref{assum:compactX}) and \ref{assum:compactOm}, and $\{Q^k\}$ is a sequence of continuous (piecewise linear and convex) functions. Further, the construction of the set of dual vertices satisfies $\Pi^0 = \{{\bf 0}\} \subseteq \ldots \subseteq \Pi^k \subseteq \Pi^{k+1} \subseteq \ldots \subseteq \Pi$ which ensures that $0 \leq Q^k(x,\obs) \leq Q^{k+1}(x,\obs) \leq Q(x,\obs)$ for all $(x,\obs)\in (\set{X}, \rvset)$. Since $\{Q^k\}$ increases monotonically and is bounded by a finite function $Q$ (due to \ref{assum:completeRecourse}), this sequence pointwise converges to some function $\xi(x,\obs) \leq Q(x,\obs)$. Once again due to \ref{assum:completeRecourse}, we know that the set of dual vertices $\Pi$ is finite and since $\Pi^k \subseteq \Pi^{k+1} \subseteq \Pi$, the set $\lim_{k\rightarrow\infty}\Pi^k := \overline{\Pi} ~(\subseteq \Pi)$ is also a finite set. Clearly, 
\begin{align*}
    \xi(x,\obs) = \lim_{k \rightarrow \infty} Q^k(x,\obs) = \max~\{\pi^\top [r(\obs) - T(\obs)x]~|~ \pi \in \overline{\Pi}\}
\end{align*}
is the optimal value of a LP, and hence, is a continuous function. The compactness of $\mathcal{X} \times \rvset$, and continuity, monotonicity and pointwise convergence of $\{Q^k\}$ to $\xi$ guarantees that the sequence uniformly converges to $\xi$ (Theorem 7.13 in \cite{Rudin1976}).
\end{proof}

\bibliographystyle{plain}
\bibliography{sampledDRO}

\end{document}